\documentclass[11pt]{article}
\usepackage{graphicx, color, amssymb}
\def\la{\big\langle}
\def\ra{\big\rangle}

\def\ds{\displaystyle}
\def\forall{\hbox{for all}~}
\def\L{{\bf L}}

\def\N{{\cal N}}
\def\I{{\cal I}}

\def\E{{\cal E}}

\def\ve{\varepsilon}

\def\G{{\cal G}}

\def\R{\mathbb R}

\def\diams{\diamondsuit}
\def\wto{\rightharpoonup}
\def\implies{\Longrightarrow}

\def\c{\centerline}

\def\vs{\vskip 2em}

\def\v{\vskip 1em}
\def\O{{\cal O}}

\def\begi{\begin{itemize}}
\def\endi{\end{itemize}}

\def\ov{\overline}
\def\Tilde{\widetilde}
\def\Hat{\widehat}
\def\bega{\begin{array}}
\def\enda{\end{array}}
\def\meas{\hbox{meas}}

\def\bel{\begin{equation}\label}
\def\eeq{\end{equation}}
\def\sqr#1#2{\vbox{\hrule height .#2pt
\hbox{\vrule width .#2pt height #1pt \kern #1pt
\vrule width .#2pt}\hrule height .#2pt }}
\def\square{\sqr74}
\def\endproof{\hphantom{MM}\hfill\llap{$\square$}\goodbreak}

\newtheorem{theorem}{Theorem}[section]
\newtheorem{lemma}[theorem]{Lemma}

\newtheorem{corollary}[theorem]{Corollary}
\newtheorem{definition}[theorem]{Definition}
\newtheorem{remark}[theorem]{Remark}

\begin{document}

\title{\bf  Weighted Irrigation Plans}
\vs

\author{Alberto Bressan and Qing Sun\\
\,
\\
Department of Mathematics, Penn State University \\
e-mails: axb62@psu.edu, ~qxs15@psu.edu}
\maketitle
\begin{abstract} 
We model an irrigation network where lower branches must be thicker in order to support the weight of the higher ones.  This leads to a countable family of ODEs,
one for each branch, that must be solved by backward induction.
Having introduced conditions that guarantee the existence and uniqueness of solutions,
our main result establishes the lower semicontinuity of the corresponding cost functional,
w.r.t.~pointwise convergence of the irrigation plans.  In turn, this yields the existence of an
optimal irrigation plan, in the presence of these additional weights.
\end{abstract}
\vs
\section{Introduction}
\label{s:1}
\setcounter{equation}{0}

In the classical irrigation problem with Gilbert cost \cite{G},
water is pumped out from a well and transported to finitely many locations $P_1,\ldots,P_n$
by a network of pipes.   The total cost is 
computed by
\bel{GC}\sum~[\hbox{flux of water through the pipe}]^\alpha\times [\hbox{length of the pipe}].\eeq
Here the sum ranges over all pipes in the network, while $\alpha\in [0, 1]$
is a fixed exponent.

This model is appropriate for an irrigation network built at ground level.
On the other hand, sometimes one would like to model a network as a free standing structure.
For example, in \cite{BS} the authors considered
tree branches transporting water and nutrients from the root to all the leaves.
In this case, one should take into account that the lower portion of each branch 
bears the weight of the upper part.
As a result, the thickness (and hence cost per unit length) of the 
lower portion should be greater, even if the flux remains the same.  
This is indeed observed in nature, where the thickness of
tree branches decreases in a continuous fashion, as one moves toward the tip.

Aim of this paper is to develop a general framework to describe this situation.
As a first step,  consider a single branch with length $\ell$,
parameterized by arc-length $s\in [0,\ell]$, oriented from the root toward the tip.
To account for the variable thickness of this branch we introduce a weight function $W=W(s)$.
Assuming that the flux is constant along the entire branch, this will satisfy
an ODE of the form
\bel{ODE} W'(s)~=~- f(W(s)),\eeq
where  $f$ is a non-negative, continuous function.
A natural set of assumptions on $f$ is
\begi
\item[{\bf (A1)}] {\it The function $f:\R_+\mapsto\R_+$ is continuous on $[0,\,+\infty[\,$, 
twice continuously differentiable for $s>0$, and satisfies
\bel{fp1} f(0)=0,\qquad f'(s)> 0,\qquad f''(s)\leq 0\qquad\forall ~s>0.\eeq
}
\endi
A typical example is $f(s) = c s^\beta$, for some $c\geq 0$ and  $0<\beta\leq 1$.
\v


\subsection{A model with finitely many branches.}

Next, we describe how to construct a family of weights $W_i(\cdot)$ in a network 
consisting of finitely many branches $\gamma_i$, $i=1,\ldots,N$.  This will be achieved by induction, starting from the tip
of each branch and proceeding backward toward the root.

\begin{figure}[ht]
\centerline{\hbox{\includegraphics[width=10cm]{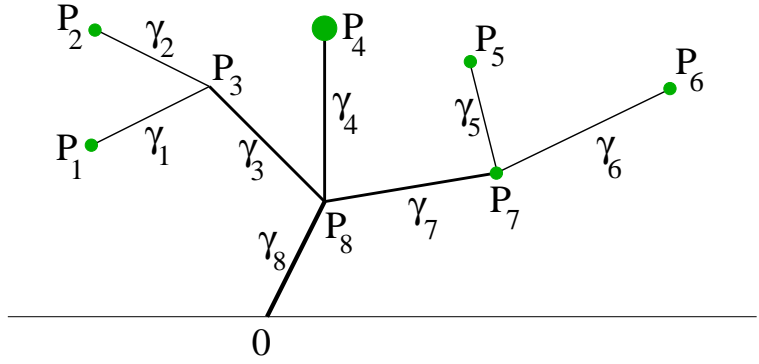}}}
\caption{\small  According to (\ref{i1n}), the branches of this tree are partitioned according to $\I_1=\{ 1,2,4,5,6\}$, $\I_2=\{3,7\}$, $\I_3=\{8\}$.
Weights can be constructed by induction, solving the backward Cauchy problems
(\ref{WPi})-(\ref{dw}) first along the arcs $\gamma_i$, $i\in \I_1$, then for $i\in \I_2$,
etc.
}
\label{f:ir92}
\end{figure}

To fix the ideas,  let each branch $\gamma_i:[0,\ell_i]\mapsto\R^d$ be parameterized by arc-length, oriented from the root toward the tip. 
As shown in Fig.~\ref{f:ir92}, call $P_i=\gamma_i(\ell_i)$ the endpoint 
of the arc $\gamma_i$ and  consider 
a measure $\mu$ consisting of finitely many point  masses $m_i\geq 0$ located 
at points $P_i$.    It is assumed that, for each node
$P_i$, there is a unique path (i.e., a concatenation of  arcs) connecting $P_i$ to the origin. 

Call 
\bel{OI}\O(i)~=~\Big\{ j\in \{1,\ldots,N\}\,;~~\gamma_j(0)= P_i\Big\}\eeq
the set of branches originating from the node $P_i=\gamma_i(\ell_i)$.
Moreover, consider the sets of indices inductively defined by
\bel{i1n}\bega{c}\I_0\,\doteq\,\emptyset,\qquad\qquad
\I_1~\doteq~\Big\{i\in \{1,\ldots,N\}\,;~~\O(i)=\emptyset\Big\},\\[3mm]
 \I_{p+1}~\doteq~\Big\{i\in \{1,\ldots,N\}\,;~~
 \O(i)\subseteq \I_1\cup\cdots\cup \I_{p}\Big\}
\setminus  (\I_1\cup\cdots\cup \I_{p})\,.\enda\eeq
Roughly speaking, $\I_1$ is the set of outer-most branches.  Branches
in $\I_p$ originate from the tips of branches in $\I_{p+1}$, etc.
Since the graph contains no loops, according to the above construction 
the set $\{1,\ldots,N\}$
is the disjoint union of the sets $\I_p$, $p\geq 1$.

For each branch $i\in \{1,\ldots,N\}$,
a weight function $W_i(\cdot)$ can now be defined in terms of the following rules:
\begi
\item[(i)]  The weight at the tip of the $i$-th branch is 
\bel{WPi}\ov W_i~\doteq~
W_i(\ell_i)~=~m_i+ \sum_{j\in \O(i)} W_j(0+).\eeq
\item[(ii)] Along each branch $\gamma_i$, the weight $W_i(\cdot)$ is absolutely continuous
and satisfies the ODE
\bel{dw}
W'_i(s)~=~-f(W_i(s))\qquad\qquad s\in \,]0, \ell_i]\,.\eeq
\endi
\v
According to (i)-(ii), the solution can be computed by induction on the entire 
tree, first on all  branches $i\in\I_1$, then on all branches $i\in \I_2$, etc.
For sake of definiteness, we  assume
\bel{m>0}m_i\,>\,0\qquad\forall i\in \I_1\,.\eeq
This guarantees that the flux along every branch is strictly positive.
In turn, by (\ref{fp1}), it implies that the backward 
Cauchy problem (\ref{WPi})-(\ref{dw}) on $[0,\ell_i]$
has a unique solution.
\v
{\bf Example 1.} When $f(s) = c s^\beta$, the Cauchy problem (\ref{WPi})-(\ref{dw}) can be solved explicitly. Namely:
$$\ov W_i^{1-\beta} - W_i^{1-\beta}(s)~=~c (1-\beta)  (s-\ell_i),$$
\bel{Wis}
W_i(s)~=~\Big( \ov W_i^{1-\beta}  +c (1-\beta)  (\ell_i-s)\Big)^{1\over 1-\beta}.\eeq
In particular, from (\ref{WPi}) we deduce the inductive rule
\bel{wind}
\ov W_i~=~m_i + \sum_{j\in \O(i)} \Big( \ov W_j^{1-\beta}  +c (1-\beta) \ell_j\Big)^{1\over 1-\beta}.
\eeq

In the presence of a weight  function $W$, for a given $\alpha\in \,]0,1]$
the {\bf total weighted cost of the irrigation network} is then defined as
\bel{EC1}
E^{W,\alpha}~\doteq~\sum_i \int_0^{\ell_i} [W_i(s)]^\alpha\, ds.\eeq
More generally, given a positive, nondecreasing, concave function 
$\psi:\R_+\mapsto\R_+$, satisfying the same assumptions as $f$ in {\bf (A1)}.
We then define
\bel{ECpsi}
E^{W,\psi}~\doteq~\sum_i \int_0^{\ell_i} \psi\bigl(W_i(s)\bigr)\, ds.\eeq
\v
\begin{remark} {\rm  In the case where $f\equiv 0$, the weight functions are constant
along every branch:
$$W_i(s)~=~\hbox{flux through}~\gamma_i\,,\qquad\qquad s\in [0, \ell_i]\,.$$
Hence the total weighted cost (\ref{EC1}) coincides with the Gilbert cost (\ref{GC}).
}\end{remark}
\v
\begin{remark} {\rm  
 In the special case where $\beta=\alpha$, so that  $f(s)= c s^\alpha$,
in view of (\ref{Wis}) this cost is computed by
\bel{EC11}
E^{W,\alpha}~\doteq~\sum_i \int_0^{\ell_i} 
\Big( \ov W_i^{1-\alpha}  +c (1-\alpha)  (\ell_i-s)\Big)^{\alpha\over 1-\alpha}\, ds.\eeq
When $\alpha=1-\alpha=1/2$, the formulas (\ref{wind}) and (\ref{EC11})
further simplify to
$$
\ov W_i~=~m_i + \sum_{j\in \O(i)} \Big( \ov W_j^{1/2}  +{c \ell_j\over 2}\Big)^2,\qquad\quad
E^{W,\alpha}~\doteq~\sum_i 
\Big( \ell_i \ov W_i^{1/2}  +{c \ell_i^2\over 4} \Big).$$
}\end{remark}

\v
Aim of the present paper is to extend the theory of optimal irrigation networks
\cite{BCM, BCM1, BS, MMS, X3, X15}, accounting for the presence of weights in the cost function.
In essence, this requires the solution of a countable family of measure-valued ODEs, one for 
each branch of the network.

In the case  of a finite network, where $\mu$ consists of finitely many atoms, 
our definition reduces to (\ref{EC1}).
For a general network, irrigating a  positive Radon measure $\mu$,
the weighted cost will be defined as a limit of an increasing sequence of 
approximations.  For any $\ve>0$, 
these approximations are obtained by restricting the transport plan to a finite set of 
paths where the flux is $\geq \ve$.

Besides showing how this family of weights can be uniquely determined,
our main results include the lower semicontinuity of the weighted irrigation cost.
As an immediate consequence, this yields the existence of optimal weighted irrigation plans.   
Furthermore, the  optimization problems for tree branches considered 
in \cite{BSun, BPS} still have solutions
when the cost functional includes the presence of weights.

The remainder of the paper is organized as follows.  Section~2 reviews  
some basic definitions and results concerning the Lagrangian approach to 
optimal irrigation plans.  
For later use, we also include some lemmas on ODEs with measure-valued
right hand side, formulated as integral  equations.
In Section~3 we provide a detailed construction of the weight functions, and the total weighted cost
of an irrigation plan.   The lower semicontinuity of the weighted cost, 
w.r.t.~pointwise convergence of the particle paths, is stated as Theorem~\ref{t:51}, and proved
in Section~4. 
In Section~5 we consider a more general model  where the increase in the  thickness of a branch, as one moves from the tip toward the root,  depends also on the inclination of the branch itself. The   
ODE (\ref{ODE})  is thus replaced by
\bel{GODE}
W'(s)~=~-f(\dot \gamma(s), W(s)),\eeq
assuming that $f$ is continuous in both variables, and that the map
$v\mapsto f(v,W)$ is positively homogeneous and convex w.r.t.~the 
variable $v\in\R^d$. We show that all previous results, including the 
lower semicontinuity of the weighted irrigation cost, remain valid in this more general case.

Finally,
in Section~6 we prove the existence of an optimal weighted irrigation plan
for a given measure $\mu$,
and the lower semicontinuity of the weighted irrigation cost w.r.t.~weak convergence of measures: $\mu_n\wto \mu$.  
In particular,  the existence of solutions 
to the optimization problems for tree branches  studied in 
\cite{BSun, BPS} remains valid also in the presence of weights.

The problem of determining which measures have a finite 
or infinite weighted irrigation cost,
depending on the dimension of their support, will be discussed in the 
forthcoming paper \cite{Sun}.  An interesting open question is whether, in the presence of weights,
an optimal irrigation plan can be computed using a suitable Modica-Mortola approximation based on 
 $\Gamma$-convergence, as in  \cite{M, OS,  PSX, S2}.

A general introduction to the theory of ramified transport can be found in \cite{BCM}.
For solutions of ODEs with measure-valued right hand side we refer to 
\cite{BR, R}.

\section{Preliminaries}
\label{s:2}
\setcounter{equation}{0}

We recall some basic definitions from \cite{BCM}.  Throughout the following, we say that a map
$\gamma(\cdot)$ is
1-Lipschitz if it is Lipschitz continuous with Lipschitz constant 1.
We denote by $\Gamma$ the set of all 1-Lipschitz maps $\gamma:\R_+\mapsto \R^d$.
By Ascoli-Arzela's theorem (see Lemma~3.4 in \cite{BCM}),
this is a compact metric space with the distance
\bel{dgg}
d(\gamma_1,\gamma_2)~\doteq~\sup_{k>1} ~{1\over k} \,\Big\{ \max_{s\in [0,k]} 
|\gamma_1(s)-\gamma_2(s)|\Big\}.\eeq
Notice that (\ref{dgg}) corresponds to the topology of uniform convergence on
compact sets.

\begin{definition}\label{d:27}
Let $\mu$ be a positive Radon measure on $\R^d$, with total mass 
$M\doteq \mu(\R^d)$.  An  {\bf irrigation plan} for $\mu$ is a function
$$\chi:[0,M] \times \R_+\mapsto \R^d,$$
measurable w.r.t.~$\xi$ and
continuous w.r.t.~$t$, with the following properties.
\begi
\item[(i)] {\bf (regularity)} For a.e.~$\xi\in [0,M]$ the map $t\mapsto  \chi(\xi,t)$ 
is 1-Lipschitz  and  eventually constant.  
Namely, there exists $\tau(\xi)\geq 0$ such that
$$\left\{\bega{cl} |\chi(\xi,t)-\chi(\xi, t')|~\leq~|t-t'|\qquad&\hbox{for all}~ t,t'\geq 0,\\[3mm]
\chi(\xi, t)~=~\chi(\xi,\tau(\xi))\qquad &\hbox{for every} ~t\geq \tau(\xi).\enda\right.$$
Throughout the following, we  denote by $\tau(\xi)$ is the smallest time $\tau$ such that
$\chi(\xi,\cdot)$ is constant for $t\geq\tau$.

\item[(ii)] {\bf ($\chi$ irrigates the measure $\mu$)} For all $\xi\in [0,M]$ one has
$\chi(\xi,0)=0$. Moreover, the push-forward of the Lebesgue measure on $[0,M]$ by the map
$\xi\mapsto \chi(\xi, \tau(\xi))$ coincides with $\mu$.    In other words,
for every Borel set $V\subseteq\R^d$ on has
\bel{PF2} 
\mu(V)~=~\meas\Big\{ \xi\in [0,M]\,;~~\chi(\xi , \tau(\xi))\in V\Big\}.\eeq
\endi
\end{definition}

\begin{remark}\label{rm1} {\rm  Relying on a theorem of Scorza-Dragoni \cite{Ku, SD}, 
one can construct a partition of the interval $[0,M]$ into countably many disjoint subsets
\bel{sdt}[0,M]~=~\left(\bigcup_{j=1}^\infty K_j\right)\cup \N,\eeq
such that
\begi
\item  each $K_j$ is compact, 
\item the set $\N$ has measure zero,
\item the restriction of $\chi$ to each product set $K_j\times \R_+$ is continuous.
\endi
}
\end{remark}
Thanks to the above construction,  measurability issues can be more easily resolved.
For example, we have
\begin{lemma}\label{l:29} Let $K\subseteq [0,M]$ be a compact subset such that $\chi$ is continuous restricted to 
$K\times \R_+$.   Then the map $\xi\mapsto \tau(\xi)$ is lower semicontinuous restricted to $K$.
\end{lemma}
\v
{\bf Proof.}
Indeed, consider a sequence $\xi_n\to \xi$ of points in $K$. If $\liminf \tau(\xi_n)=+\infty$ 
there is nothing to prove. Otherwise,
by taking a subsequence, we can assume
$$\lim_{n\to\infty}\tau(\xi_n)~=~\bar\tau.$$
By assumption, the  continuous functions $\chi(\xi_n,\cdot)$ converge to $\chi(\xi,\cdot)$
uniformly on compact sets.   For any $\ve>0$, all but finitely many of these functions 
are constant on $[\bar\tau+\ve, \,+\infty[\,$.  Hence also $\chi(\xi,\cdot)$ is constant on this same 
domain.   Since $\ve>0$ is arbitrary, we conclude that $\chi(\xi,\cdot)$ is constant 
on $[\bar\tau,+\infty[\,$.   Hence $\tau(\xi)\leq\bar\tau$, as claimed.
\endproof
\v
\begin{corollary}\label{c:210}
 Given any $\ve>0$ there exists a compact set $K\subseteq [0,M]$, with
\bel{MK}
\meas\Big([0,M]\setminus K\Big)~<~\ve,\eeq
and such that 
\begi
\item[(i)] the map $\xi\mapsto \chi(\xi,\cdot)$ is continuous restricted to $K$,
\item[(ii)] the map $\xi\mapsto\tau(\xi) $ is continuous restricted to $K$.
\endi
\end{corollary}

{\bf Proof.} By Remark~\ref{rm1},  we can choose $K=\cup_{j=1}^\nu K_j$ with $\nu$ 
large enough so that (\ref{MK}) holds. Since $\chi$ is continuous on each $K_j$, 
the statement (i) follows immediately.   By 
Lemma~\ref{l:29} the map $\xi\mapsto \tau(\xi)$ is measurable on $K$.
By Lusin's theorem, there exists a smaller compact set $K_0\subseteq K$,
still with meas$([0,M]\setminus K_0)~<~\ve$,
such that the restriction of $\tau(\cdot)$ to $K_0$ is continuous.
By replacing $K$ with $K_0$, the conclusion (ii) of the Corollary is satisfied.
\endproof
\v

The usual definition of {\em irrigation cost} 
involves the multiplicity of a point $x=\chi(\xi,t)$, defined as
\bel{chim}
|x|_\chi~\doteq~\meas\Big(\{\xi\in [0,M]\,;~~\chi(\xi,t)=x~~\hbox{for some}~t\geq 0\}\Big). \eeq
In the present case, this must be replaced by a different concept, 
related to the single-path
property.

\begin{definition}  
\label{d:eq1}
We say that two 1-Lipschitz maps $\gamma:[0,t]\mapsto \R^d$ and $\gamma':[0, t']\mapsto\R^d$ 
are {\bf equivalent} if they are parameterizations of the same curve.  That is, if
there exists an interval $[0,T]$ and nondecreasing, Lipschitz continuous surjective maps
$s\mapsto \eta(s)$ and $s\mapsto \eta'(s)$ from $[0,T]$ onto $ [0,t]$ and $[0, t']$ respectively, such that
\bel{eq1}
\gamma(\eta(s))~=~\gamma'(\eta'(s))\qquad\qquad\forall s\in [0,T].\eeq
If this is the case, we write $\gamma\simeq\gamma'$.
\end{definition}

\begin{remark}  {\rm Given a 1-Lipschitz map $\gamma:[0,t]\mapsto \R^d$, its arc-length
re-parameterization is the map 
$$\sigma~\mapsto~\gamma(s(\sigma))$$ where, for every $\sigma$,
one has
$$\int_0^{s(\sigma)}  |\dot \gamma(\zeta)|\, d\zeta ~=~\sigma.$$
According to the above definition,
two maps $\gamma:[0,t]\mapsto \R^d$ and $\gamma':[0, t']\mapsto\R^d$
are equivalent  if and only if their arc-length re-parameterizations coincide.
}
\end{remark}

\begin{remark}\label{r:1lip}{\rm
In  Definition~\ref{d:eq1}, one can always take $T=t+t'$ and 
assume that both functions $\eta,\eta'$ are 1-Lipschitz. 
Indeed, let $\eta,\eta'$ be maps from $[0,T]$
	onto $[0,t]$ and $[0,t']$ respectively, such that (\ref{eq1}) holds.
	
	For $s\in [0,T]$, define the nondecreasing, surjective map 
	$\sigma:[0,T]\mapsto 
	[0, \Tilde T]\doteq [0, t+t']$ by setting
	$$\sigma(s)~\doteq~\eta(s)+\eta'(s).$$
We then define the maps $\tilde\eta,\tilde\eta'$ from $[0,\Tilde T]$ into $[0,t]$ and $[0,t']$
implicitly, by setting
\bel{teta}\tilde\eta(\sigma(s))~\doteq~\eta(s),\qquad\quad \tilde\eta'(\sigma(s))~\doteq~\eta'(s)
\qquad s\in [0,T].\eeq
We claim that $\tilde\eta$ and $\tilde \eta'$ are 1-Lipschitz.  Indeed,
let $\sigma_1 =\sigma(s_1)<\sigma(s_2)=\sigma_2$.   Then
$$\bigl[\tilde \eta(\sigma_2)-\tilde \eta(\sigma_1)\bigr] + [\tilde \eta(\sigma_2)-\tilde \eta(\sigma_1)]
~=~[\eta(s_2)+\eta'(s_2)] -[\eta(s_1)+\eta'(s_1)] ~=~
\sigma_2-\sigma_1\,.$$
The identity 
$\gamma(\tilde\eta(s))~=~\gamma'(\tilde\eta'(s))$
 now follows from (\ref{eq1}) and (\ref{teta}).}
\end{remark} 
Throughout the following, 
we denote by $\gamma\Big|_{[0,t]}$ the restriction of a map $\gamma$ to the interval
$[0,t]$.
	
\begin{definition}\label{d:1}
Let $\chi:[0,M]\times \R_+\mapsto \R^d$ be an irrigation plan for the measure $\mu$.
We define an equivalence relation on the set $[0,M]\times \R_+$ by setting $(\xi,t)\,\sim\,(\xi', t')$
whenever $\chi(\xi,\cdot)\Big|_{[0,t]}\simeq\chi(\xi',\cdot)\Big|_{[0,t']}$. This means that
the maps
$$s~\mapsto~\chi(\xi,s), \quad s\in [0,t]\qquad\hbox{and}
\qquad s~\mapsto~\chi(\xi',s), \quad s\in [0,t']
$$
are equivalent in the sense of Definition~\ref{d:eq1}.

The {\bf multiplicity} of $(\xi,t)$ is then defined as
\bel{mtx}m(\xi,t)~\doteq~\meas\Big(\bigl\{ \xi'\in[0,M]\,;~~(\xi', t')\sim (\xi,t) ~~\hbox{for some}~t'>0\bigr\}
\Big).\eeq
\end{definition}

\begin{remark} {\rm The multiplicity  $m(\xi,t)$ measures the total amount of particles that 
pass through the point $x=\chi(\xi,t)$ {\em traveling along exactly the same path
as the particle $\xi$.}    If $\chi$ has the single path property 
(see Chapter 7 in \cite{BCM}), then
$m(\xi,t)=|\chi(\xi,t)|_\chi$.   However, for a general irrigation plan
we only have the inequality
\bel{mineq}
m(\xi,t)~\leq~|\chi(\xi,t)|_\chi\,.\eeq   Notice that one may well have
$$\chi(\xi,t) \,=\, \chi(\xi', t')\qquad\hbox{but}\qquad m(\xi,t)\,\not=\, m(\xi', t').$$}
\end{remark}
Given an irrigation plan $\chi:[0,M]\times\R_+\mapsto\R^d$, throughout the following
we shall use a basic assumption,
which is needed in order that the total cost be finite. As in part (i) of Definition~\ref{d:27},
we denote by $\tau(\xi)$ the time when the particle $\xi$ reaches its final 
destination.
\begi
\item[{\bf (A2)}] {\it
For a.e.~$\xi\in [0,M]$,  one has $m(\xi,t)>0$
for every $0\leq t<\tau(\xi)$. }
\endi
In other words, for any  particle $\xi$ and 
 any $t\in [0, \,\tau(\xi)[\,$, there is a positive amount
of other particles that travel along the same path $\chi(\xi,\cdot)\Big|_{[0,t]}$.

The next two lemmas establish various properties of  the multiplicity function 
introduced in Definition~\ref{d:1}

\begin{lemma} \label{l:2} For any $\ve>0$ there exists a compact set $K\subseteq [0,M]$
satisfying (\ref{MK}), such that  the set-valued function 
\bel{Fdef}F(\xi,t)~\doteq~\{\xi'\in K\,;~~(\xi', t')\sim(\xi,t)~~
\hbox{for some}~~t'\geq 0\}\eeq
is upper semicontinuous on $K\times \R_+$.
\end{lemma}

{\bf Proof.}   {\bf 1.} Given $\ve>0$, let $K\subseteq [0,M]$ be the compact set constructed
in Corollary~\ref{c:210}.
 We claim that the graph of $F$, restricted to $K\times \R_+$, is closed.
In other words, assume that 
$$\xi_n\to \xi\,\qquad  t_n\to t,\qquad \xi'_n\to \xi' \qquad\hbox{as}\quad n\to \infty,$$
and moreover 
$$(\xi'_n, t'_n)\sim(\xi_n, t_n)\qquad\forall n\geq 1.$$
We need to show that there exists $t'\geq 0$ such that $(\xi, t)\sim(\xi', t')$.
\v
{\bf 2.} By the assumptions, according to Remark~\ref{r:1lip} for every $n\geq 1$ 
there exists an interval $[0, T_n] = [0, t_n+t'_n]$ and two
 nondecreasing, 1-Lipschitz, surjective  maps
$$\eta_n: [0, T_n]\,\mapsto [0, t_n],\qquad \eta'_n: [0, T_n]\,\mapsto [0, t'_n],
$$
such that 
\bel{eat0}\chi(\xi_n, \eta_n(s))~=~\chi(\xi'_n, \eta'_n(s))\qquad\forall s\in [0, T_n].\eeq
\v
{\bf 3.} We now observe that, since the map $\xi\mapsto \tau(\xi)$ is continuous on the compact set $K$, it is uniformly bounded.   We can thus assume that the sequence 
$(t'_n)_{n\geq 1}$ is bounded.
Since $t_n\to t< +\infty$ and $T_n=t_n + t_n'$, we have the uniform boundedness of 
$(T_n)_{n\geq 1}$. By extracting a subsequence, one can assume 
\bel{eat1}	\lim_{n\to +\infty} t_n'~=~t',\qquad\lim_{n\to +\infty} T_n ~=~T\, .\eeq
If $T_n< T$, we extend the maps $\eta_n$ and $\eta_n'$ to the interval $[0, T]$ by setting 
$\eta_n(s)\doteq \eta_n(T_n)$, $\eta_n'(s)\doteq \eta_n'(T_n)$, for all $s\in (T_n, T]$. 
Using the Ascoli-Arzel\`a theorem, by possibly extracting a subsequence 
we achieve the uniform convergence
\bel{eat3}\eta_n(\cdot)\to \eta(\cdot),\qquad \eta_n'(\cdot)\to \eta'(\cdot)\, ,
\qquad\hbox{ uniformly on } ~[0, T]\, .\eeq
Here $\eta$ and $\eta'$ are two nondecreasing, 1-Lipschitz, surjective maps from 
$[0, T]$ onto $[0,t]$ and $[0, t']$ respectively. By the continuity of $\chi(\cdot,\cdot)$ 
on $K\times \R_+$, from (\ref{eat0}) we obtain
\bel{eat4}\chi(\xi, \eta(s))~=~\lim_{n\to \infty} \chi(\xi_n, \eta_n(s))~=~\lim_{n\to \infty} 
\chi(\xi_n', \eta_n'(s))~=~\chi(\xi', \eta'(s))\qquad \hbox{ for all } s\in [0, T].\eeq
Therefore, $(\xi', t')\sim (\xi, t)$.
\endproof

\begin{lemma}\label{l:31}
Let $\chi:[0,M]\times\R_+\mapsto \R^d$ be an irrigation plan for the measure $\mu$.
Then the following holds.
\begi
\item[(i)] The map $(\xi,t)\mapsto m(\xi,t)$ is measurable.
\item[(ii)]
For each $\xi\in [0,M]$,  the map $t\mapsto m(\xi,t)$  is  non-increasing and  left continuous.
\item[(iii)]
For any fixed $\ve>0$,  the stopping time
\bel{ted}
\tau_\ve(\xi)~\doteq~\max\,\Big\{t\in [0,  \tau(\xi)]\,;~~m(\xi,t)\geq \ve\Big\}\eeq
is a measurable function of $\xi\in [0,M]$.
\endi
\end{lemma}

{\bf Proof. 1.} Given $\ve>0$, let $K\subseteq [0,M]$ be a compact set 
satisfying the conditions in Lemma~\ref{l:2}.   In terms of the multifunction $(\xi,t)\mapsto F(\xi,t)\subseteq K$ defined
at (\ref{Fdef}),  this implies  the scalar
function
\bel{peo}(\xi,t)~\mapsto~\meas(F(\xi,t))\eeq
is upper semicontinuous restricted to  $K\times \R_+\,$. 
For $(\xi,t)\in K\times \R_+$  this implies
\bel{mix}m(\xi,t)-\ve~\leq~ \meas(F(\xi,t))~\leq~m(\xi,t).\eeq
\v
{\bf 2.} Repeating the above construction for decreasing values of $\ve$,
we can find an increasing sequence of compact sets $(K_n)_{n\geq 1}$, 
with $\meas([0,M]\setminus K_n)<1/n$, such that 
\bel{mm}m(\xi,t) - {1\over n}~\leq~\meas(F_n(\xi,t))~\leq~m(\xi,t)\, .\eeq
Here $F_n$ is the multifunction defined at  (\ref{Fdef}), with $K$ replaced by $K_n$.
Notice that the function $(\xi,t)~\mapsto~\meas(F_n(\xi,t))$
is upper semicontinuous restricted to     $K_n\times \R_+$.  
Setting
$$m_n(\xi,t)~\doteq~\left\{ \bega{cl}\meas(F_n(\xi,t))&\qquad\hbox{if}
~~\xi\in K_n\,,\\[3mm]
0&\qquad\hbox{if}
~~\xi\notin K_n\,,\enda\right.$$
by (\ref{mm}) we have the pointwise convergence $m_n(\xi,t)\to m(\xi,t)$
for a.e.~$\xi\in [0,M]$.   Since each $m_n$ is measurable, the same holds for $m$.
This proves (i).
\v
{\bf 3.} By the definition of the multiplicity function in (\ref{mtx}), it immediately follows that the map $t\mapsto m(\xi,t)$ is non-increasing. To prove its left continuity, fix $(\xi,t)\in [0,M]\times \R_+$ and consider an increasing  sequence $t_n\uparrow t$.
By monotonicity, it  follows 
 \bel{co00}\lim_{n\to\infty}m(\xi,t_n)~=~\inf_n \, m(\xi, t_n)~\geq~m(\xi,t)\, .\eeq
 To prove that equality actually holds in (\ref{co00}), 
given any $\ve>0$, let $K\subseteq [0,M]$ be a compact set satisfying the conditions in Lemma \ref{l:2}.   By the upper semicontinuity of the multifunction  $t\mapsto F(\xi,t)$
one has
\bel{co12}m(\xi,t)~\geq~\meas(F(\xi,t))~\geq~\limsup_{n\to\infty}~\meas(F(\xi,t_n))~\geq~\limsup_{n\to\infty}~m(\xi,t_n) - \ve\, . \eeq
Since $\ve>0$ was arbitrary, this   proves statement (ii) of the lemma.
\v
{\bf 4.} To prove (iii), we first observe that, by the arguments in the previous
steps {\bf 1 - 2}, for each fixed $t>0$ the map 
$\xi\mapsto m(\xi,t)$ is measurable. 
Moreover, by
 Corollary \ref{c:210} it follows that $\xi\mapsto \tau(\xi)$ is measurable. 
 For every $t>0$ we have the identity 
\bel{stpn}\Big\{\xi\in [0,M]\, ;~\tau_\ve(\xi) \geq t\Big\}~=~\Big\{\xi\in [0,M]\, ;~m(\xi,t) \geq \ve\Big\} \bigcap \Big\{\xi\in [0,M]\, ;~\tau(\xi) \geq t\Big\}.\eeq
This implies that the map $\tau_\ve(\cdot)$ is measurable.
\endproof
\v
\subsection{ODE's with measure-valued right hand side.}

For future use, we now prove some results on existence and continuous dependence, 
for Carath\'eodory solutions to an ODE  backward in time. Since in our equations
the  right hand side can 
possibly be a measure, it will be convenient to study  directly the corresponding
 integral equations.

\begin{lemma}\label{l:mvode}
Let $f:[\ve, +\infty[\,\mapsto\R_+$ be Lipschitz continuous. 
For $t\in [0,T]$,  let $t\mapsto m(t)$
be a non-increasing function with $m(T)\geq \ve$.  

\begi
\item[(i)] There exists a unique function
$w:[0,T]\mapsto [\ve, +\infty[\,$ which satisfies the integral equation
\bel{ie3}
w(t)~=~\int_t^T f(w(s))\, ds + m(t)\qquad\qquad \forall t\in [0,T].\eeq
\item[(ii)] If $m_1(t)\leq m_2(t)$ for all $t\in [0,T]$, then the corresponding
solutions of (\ref{ie3}) satisfy
\bel{w12}
w_1(t)~\leq~w_2(t)\qquad\qquad\forall t\in [0,T].\eeq
\item[(iii)] Consider a sequence of measurable sets $J_n\subseteq [0,T]$ 
such that $\lim_{n\to\infty}\meas(J_n) = 0$, and define the functions
$$ f_n( t, \omega )~\doteq~\left\{\bega{cr} f(\omega)&\qquad \hbox{if}~~t\notin J_n\,,\\[3mm]
0&\qquad \hbox{if}~~t\in J_n\,. \enda\right.  $$
Let
$t\mapsto m_n(t)\in [\ve,+\infty[\,$ be a sequence of non-increasing functions such that,
as $n\to\infty$, 
\bel{mmn} \|m_n-m\|_{\L^1([0,T])}\to 0.
\qquad 
\qquad
m_n(0)\to m(0+).\eeq
Then the solutions to
 \bel{ie30}
w_n(t)~=~
\int_t^T f_n(s,\omega_n(s))\, ds + m_n(t)\qquad\qquad \forall t\in [0,T].\eeq
satisfy
\bel{wwn} \|w_n-w\|_{\L^1([0,T])}\to 0,\qquad 
\qquad
w_n(0)\to w(0+).\eeq
\endi
\end{lemma}

{\bf Proof.}  {\bf 1.} Consider the function
\bel{FFn}F(t,z)\,\doteq\, f(m(t)+z).\eeq
We observe that a map $t\mapsto w(t)$ satisfies the integral 
equation (\ref{ie3}) if and only if
$z(t)= w(t)-m(t)$ provides a Carath\'eodory solution to the backward Cauchy problem
\bel{CPT}
\dot z(t)~=~-F(t, z(t)),\qquad\qquad z(T)=0.\eeq
Observing that $F$ is measurable in $t$ and uniformly Lipschitz continuous
in $z$, by the standard theory of  ODE  \cite{H} we conclude that (\ref{CPT})
has a unique  solution $t\mapsto z(t)$.
In turn, $w(t) = z(t)+ m(t)$ provides the unique solution to (\ref{ie3}).
\v
{\bf 2.} To prove (ii), for $i=1,2$ let $z_i$ be a solution to
$$-\dot z_i(t)~=~F_i(t, z_i(t))~\doteq~f(m_i(t)+z_i(t)),\qquad\qquad z_i(T)=0.$$
Since $F_1(t,z)\leq F_2(t,z)$ for all $t,z$, and both $F_1,F_2$ are Lipschitz
continuous w.r.t.~$z$, a standard comparison argument yields
$z_1(t)\leq z_2(t)$ for all $t\in [0,T]$.   In turn this implies
$$w_1(t)~=~m_1(t) +z_1(t)~\leq~m_2(t) +z_2(t)~=~w_2(t).$$
\v
{\bf 3.} To prove (iii), set $ F_n(t,z)\doteq f_n(t,m(t) + z)$
and let $z_n$ be the solution to 
\bel{CPn}
\dot z_n(t)~=~-F_n(t, z_n(t)),\qquad\qquad z_n(T)=0.\eeq
Then the difference $\eta_n(t)\doteq |z_n(t)-z(t)|$ satisfies
$$\bega{l}\ds\eta_n(t)~\leq~\int_t^T \Big| f_n(s,m_n(s) + z_n(s)) - f(m(s) + z(s) )  \Big|\, ds\\[3mm]
\ds~\leq~\int_t^T\Big|f(m_n(s)+ z_n(s)) - f(m(s)+ z(s)) \Big|  +  \chi_{J_n}(s)\cdot\Big|f(m(s)+ z(s))\Big|\, ds\\[3mm]
\ds~\leq~\int^T_t L\eta_n(s) + L|m_n(s) - m(s)| + \chi_{J_n}(s)\cdot|f(m(s)+ z(s)|\, ds\,.\enda$$
Here $L$ is a Lipschitz constant for the function $f$ on $[\ve,+\infty[\,$, while
$\chi_{J_n}$ denothes the characteristic function of the set $J_n$.
By Gronwall's inequality one obtains
\bel{ete}
\eta_n(t)~\leq~\int_t^T e^{L(s-t)}\Big[ L |m_n(s)-m(s)| + \chi_{J_n}(s)\cdot|f(m(s) + z(s) )|\Big]\,ds\,.\eeq
Since multiplicity functions are non-increasing, there exists some finite constant $K>0$ such that
\bel{ete0} |f(m(s) +z(s) )|~\leq~K\, ,\qquad \forall n\geq1,~s\in [0,\, T]\, . \eeq
Letting $n\to\infty$, by (\ref{mmn}) and (\ref{ete})-(\ref{ete0}), since $\lim_{n\to \infty} \meas(J_n )=0$,  we thus have the convergence
$\eta_n(t)=|z_n(t)-z(t)|\to 0$  uniformly for $t\in [0,T]$.
Recalling that $w_n = z_n + m_n$ and $w=z+m$, from (\ref{mmn}) 
it now follows (\ref{wwn}).
\endproof

	\begin{lemma}\label{l:simple}
		Let $m, m_1,\ldots, m_q:[0,\ell]\mapsto [\ve_0, +\infty[\, $ be 
non-increasing  functions such that 
\bel{nk1}\sum_{i=1}^q m_i(s)~\geq~m(s),\qquad\qquad \forall s\in [0,\ell].\eeq
Assume that  $f$ satisfies {\bf (A1)} and let $w,w_i :[0,\ell]\mapsto [\ve_0, +\infty[$,
be solutions to 
\bel{nk3}w(s)~=~\int^{\ell}_s f(w(t))\, dt + 
m(s),\qquad w_i(s)~=~\int^\ell_s f(w_i(t))\, dt + m_i(s),\eeq	
respectively.   Then,  for all  $s\in [0,\ell]$, one has
\bel{nk2}\sum_{i = 1}^q w_i(s)~\geq~w(s) .\eeq
	\end{lemma}

	{\bf Proof. } Consider the functions
$$\Tilde w(s)~\doteq~\sum_{i=1}^q w_i(s),\qquad\Tilde m(s)~\doteq~\sum_{i=1}^q
m_i(s).$$
Using the properties (\ref{fp1})  of the function $f$ and the inequality (\ref{nk1}),
 for all $s\in [0,\ell]$,
$$\Tilde w(s)~=~\int_s^\ell \sum_i f(w_i(t)\, dt + \sum_i m_i(s)
~\geq~\int_s^\ell  f(\Tilde w(t))\, dt +   \Tilde m(s).$$
Since $\Tilde m(s)\geq m(s)$, the comparison property stated in (iii) of Lemma~\ref{l:mvode} now implies $\Tilde w(s)\geq w(s)$ for all $s\in [0,\ell]$.
\endproof

\section{Construction of the weight functions}
\label{s:23}
\setcounter{equation}{0}
Given an irrigation plan $\chi:[0,M]\times\R_+\mapsto\R^d$ and a function $f$ satisfying
{\bf (A1)}, 
in this section we construct the weight function $W=W(\xi,t)$, 
by taking the supremum of a family
of approximations $W^\ve$.   

Recalling the equivalence relation introduced in  Definition~\ref{d:eq1},
we introduce

\begin{definition}\label{d:good}
Given an irrigation plan $\chi$, we say that a path $\gamma:[0, \ell]\mapsto\R^d$,
parameterized by arc-length,
is {\bf $\ve$-good} if 
\bel{goode}\meas\left(\Big\{\xi\in[0,M]\,;~~\chi(\xi,\cdot)\Big|_{[0,t]} \simeq
~\gamma ~~\hbox{for some}~ t= t(\xi)>0\Big\}
\right)~\geq~\ve.\eeq
The family of all $\ve$-good paths will be denoted by $\G_\ve$.
\end{definition}
In other words,  $\gamma$ is $\ve$-good if there is an amount  $\geq \ve$
of particles whose trajectory contains $\gamma$ as its initial 
portion. A somewhat similar definition can be found in \cite{S1}.

The family of all curves parameterized by arc-length comes with a natural partial order.
Namely, given two maps $\gamma:[0,\ell]\mapsto\R^d$,
$\gamma':[0,\ell']\mapsto\R^d$,  
we write $\gamma\preceq \gamma'$  if $\ell\leq \ell'$ and $\gamma'(s)=\gamma(s)$ for
all $s\in [0,\ell]$.   
The next lemma yields a bound on the number of maximal curves, within the family
of $\ve$-good paths.
\v
\begin{lemma} \label{l:210}
Given an irrigation plan $\chi:[0,M]\times\R_+\mapsto \R^d$ and  $\ve>0$,
 there can be at most $M/\ve$ distinct maximal $\ve$-good paths.\end{lemma}

{\bf Proof.} Let $\gamma_1,\ldots,\gamma_\nu$ be distinct maximal $\ve$-good paths.
For each $i\in\{1,\ldots,\nu\}$, consider the set
\bel{Ai}A_i~\doteq~ \Big\{\xi\in[0,M]\,;~~\chi(\xi,\cdot)\Big|_{[0,t]} \simeq
~\gamma_i ~~\hbox{for some}~  t>0\Big\}.\eeq
We claim that all sets $A_i$ are disjoint.   Indeed, if $\xi\in A_i\cap A_j$, then 
$$\chi(\xi,\cdot)\Big|_{[0,t]}\simeq\gamma_i\,,\qquad\chi(\xi,\cdot)\Big|_{[0,t']}\simeq\gamma_j\,.
$$
To fix the ideas, assume $t\leq t'$.    Then $\gamma_i\prec\gamma_j$, against the maximality 
of $\gamma_i$. This contradiction proves our claim.   In turn this implies $\nu\leq M/ \ve$,
proving the lemma.
\endproof
\v
We now fix $\ve>0$, and let  $\{\Hat \gamma_1,\ldots,\Hat \gamma_\nu\}$ be the set of all 
maximal $\ve$-good paths for the irrigation plan $\chi$.
Along each path $\Hat \gamma_i:[0,\hat \ell_i]\mapsto\R^d$ we define the {\bf multiplicity} $\Hat m_i: [0, \hat \ell_i]
\mapsto \R_+$
by setting 
\bel{mi}
\Hat m_i(t)~\doteq~\meas\bigg( \left\{ \xi\in [0,M]\,;~~
\hbox{there exists $t'\geq 0$ such that }~\chi(\xi,\cdot)\Big|_{[0,t']}~\simeq~
\Hat \gamma_i\Big|_{[0,t]}\right\}\bigg).\eeq
Otherwise stated, $\Hat m_i(t)$ is the amount of particles that travel
along the path $\Hat \gamma_i$, at least up to the point $\Hat \gamma_i(t)$. 

To construct the weight functions, we first need to split the 
maximal paths $\Hat\gamma_i$ into elementary paths $\gamma_k$, to which
an inductive procedure as in (\ref{WPi})-(\ref{dw})
can then be applied.
With this goal in mind, we define the bifurcation times
\bel{tauij}
\tau_{ij}~=~\tau_{ji}~\doteq~\max\,\Big\{ t\geq 0\,;~~\Hat\gamma_i(s)=\Hat\gamma_j(s)~~\forall
s\in [0,t]\Big\}.\eeq
The elementary paths $\gamma_k:[a_k, b_k]\mapsto \R^d$ and the corresponding 
multiplicity 
functions $m_k$ are constructed  by the following Path Splitting Algorithm.
\begi
\item[{\bf (PSA)}] For each $i\in\{1,\ldots,\nu\}$,  consider the set 
$$\{\tau_{i1}, \ldots,\tau_{i\nu}\} ~=~\{t_{i,1},\ldots, t_{i,N(i)}\},$$
where the times
\bel{tij}0\,< \,t_{i,1}\,<\,t_{i,2}\,<\,\cdots\,<\,t_{i,N(i)}\, =\,\hat \ell_j\eeq
provide an increasing arrangement of the set of  times $\tau_{ij}$ where the path
$\Hat\gamma_i$ splits apart from other maximal paths.
For each $k=1,\ldots, N(i)$, let $\gamma_{i,k}$
be the restriction of the maximal path $\Hat \gamma_i$ to the subinterval $[t_{i,k-1}, t_{i,k}]$.
The multiplicity function $m_{i,k}$ along this path is defined simply as
\bel{multj}
m_{i,k}(t)~=~\Hat m_i(t)\qquad\qquad t\in [t_{i,k-1}, t_{i,k}].\eeq

If $\tau_{ij}>0$, i.e.~if the two maximal paths $\Hat \gamma_i$ and 
$\Hat \gamma_j$ partially overlap, it is clear that  some of the elementary  paths $\gamma_{i,k}$
will coincide with some $\gamma_{j,l}$.
To avoid listing multiple times the same path, we thus remove from our list 
all paths
$\gamma_{j,l}:[t_{j, l-1}, t_{j,l}]\mapsto\R^d$ such that 
$t_{j,l} \leq \tau_{ij}$ for some $i<j$.    
After relabeling all the remaining paths, the algorithm yields a 
family of elementary paths and corresponding multiplicities 
\bel{eph}\gamma_i:[a_i,b_i]\mapsto \R^d,\qquad m_i:\, [a_i,b_i]\mapsto \R_+\,,
\qquad i=1,\ldots,N.\eeq
 \endi

For example, the tree shown in Fig.~\ref{f:ir92} contains 5 maximal paths
$\Hat\gamma_1,\ldots, \Hat\gamma_5$.   These can be
decomposed it into 8 elementary paths $\gamma_1,\ldots,\gamma_8$.
Each maximal path is a concatenation of elementary paths, namely
$$\Hat\gamma_1= \gamma_8\circ\gamma_3\circ \gamma_1\,,\qquad 
\Hat\gamma_2= \gamma_8\circ\gamma_3\circ \gamma_2\,,\qquad\Hat\gamma_3= \gamma_8\circ\gamma_4 \,,\quad \ldots $$

A set of weight functions $W_i$ on the elementary branches $\gamma_i$ can 
now be constructed
by a backward inductive procedure, similar to (\ref{WPi})-(\ref{dw}).
As in (\ref{OI}), call
$\O(i)$
the set of branches originating from the node $P_i=\gamma_i(b_i)$.
Moreover, consider the sets of indices $\I_p$ inductively defined at (\ref{i1n}).
\begi
\item[(i)] For $p=1$,  on each elementary path $\gamma_i:[a_i, b_i]\mapsto\R^d$ with $i\in \I_1$, the weight
$W_i^\ve(t)$ is defined to be the solution of
\bel{W1}
w(t)~=~\int_t^{b_i} f(w(s))\, ds + m_i(t),\qquad\qquad t\in \,]a_i, b_i].
\eeq
\item[(ii)]
Next, assume that
the weight functions $W^\ve_k(t)$ have already been constructed
along all paths $\gamma_k:[a_k, b_k]\mapsto\R^d$ with $k\in \I_{p-1}$.   
 
 For  $i\in \I_p$, the weight $W^\ve_i(t)$ along the 
 $i$-th branch is then defined to be the solution of
 \bel{Wilast}
w(t)~=~\int_{t}^{b_i} f(w(s))\, ds + m_i(t) + \ov w_{i}\,,
\qquad\qquad t\in \,]a_i, b_i].
\eeq
where
\bel{posti}\ov w_{i}~\doteq~\sum_{k\in \O(i)} W^\ve_k(a_k+) -\sum_{k\in \O(i)} m_k(a_k+).
\eeq
\endi
Notice that  {\bf (PSA}) implies
$b_i=a_k$ for all $k\in \O(i)$. At the end-point $\gamma_i(b_i)$,  
the weight is
$$W^\ve_i(b_i)~=~\sum_{k\in \O(i)} W^\ve_k(a_k+) +\left[m_i(b_i)
-\sum_{k\in \O(i)} m_k(a_k+)\right].$$
Here the term between brackets can be strictly positive. For example, this 
will happen
 if the irrigated measure
$\mu$ contains a point mass at $\gamma_i(b_i)$.

By induction on $p$, after finitely many steps we  obtain a weight function 
$W_i^\ve:[a_i, b_i]\mapsto [\ve, +\infty[\,$ defined on each elementary  path
$\gamma_i$.   

Going back to the maximal paths $\Hat \gamma_j$ considered 
in {\bf (PSA)}, 
the above construction yields a weight $\Hat W_{j,k}$ on the restriction 
of  $\Hat \gamma_j$ to each subinterval $[t_{j, k-1}, t_{j,k}]$.
 Along the maximal path 
$\Hat\gamma_j$, the weight $\Hat W_j:[0, \hat\ell_j]\mapsto \R_+$
is  then defined simply by setting
\bel{HWI}\Hat W_j(t)~=~\Hat W_{j,k}(t)\qquad\qquad \hbox{if}~~~
t\in [t_{j, k-1}, t_{j,k}].\eeq

\v
Next, 
in order to construct  an approximate 
weight function $W^\ve:[0,M]\times
\R_+\mapsto \R_+$ on the family of all paths $\chi(\xi,\cdot)$ of the irrigation 
plan, we consider the {\bf stopping time}
\bel{tep}
\tau_\ve(\xi)~=~\sup\bigl\{ t\geq 0\,;~~m(\xi,t)\geq\ve
\bigr\}.\eeq
We then define 
the corresponding weight function 
\bel{We1}
W^\ve(\xi,t)~\doteq~\left\{\bega{cl}\Hat W_i(s)\qquad &\hbox{if}\qquad t\leq \tau_\ve(\xi),\qquad 
\chi(\xi,\cdot)\Big|_{[0,t]}
~\simeq~\Hat \gamma_i\Big|_{[0,s]}\,,\\[4mm]
0\qquad &\hbox{if}\qquad t>\tau_\ve(\xi).\enda\right.\eeq
Having constructed these approximate weights $W^\ve$,
the weight function $W$ is then obtained by letting $\ve\to 0$. 

\begin{definition} Let $\chi:[0,M]\times \R_+\mapsto\R^d$ be an irrigation plan 
satisfying {\bf (A2)}.
The {\bf weight function} $W=W(\xi,t)$ for $\chi$ is  defined as
\bel{W}
W(\xi,t)~\doteq~\sup_{\ve >0} ~W^\ve(\xi,t).\eeq
\end{definition}

\begin{figure}[ht]
\c{\includegraphics[scale=0.6]{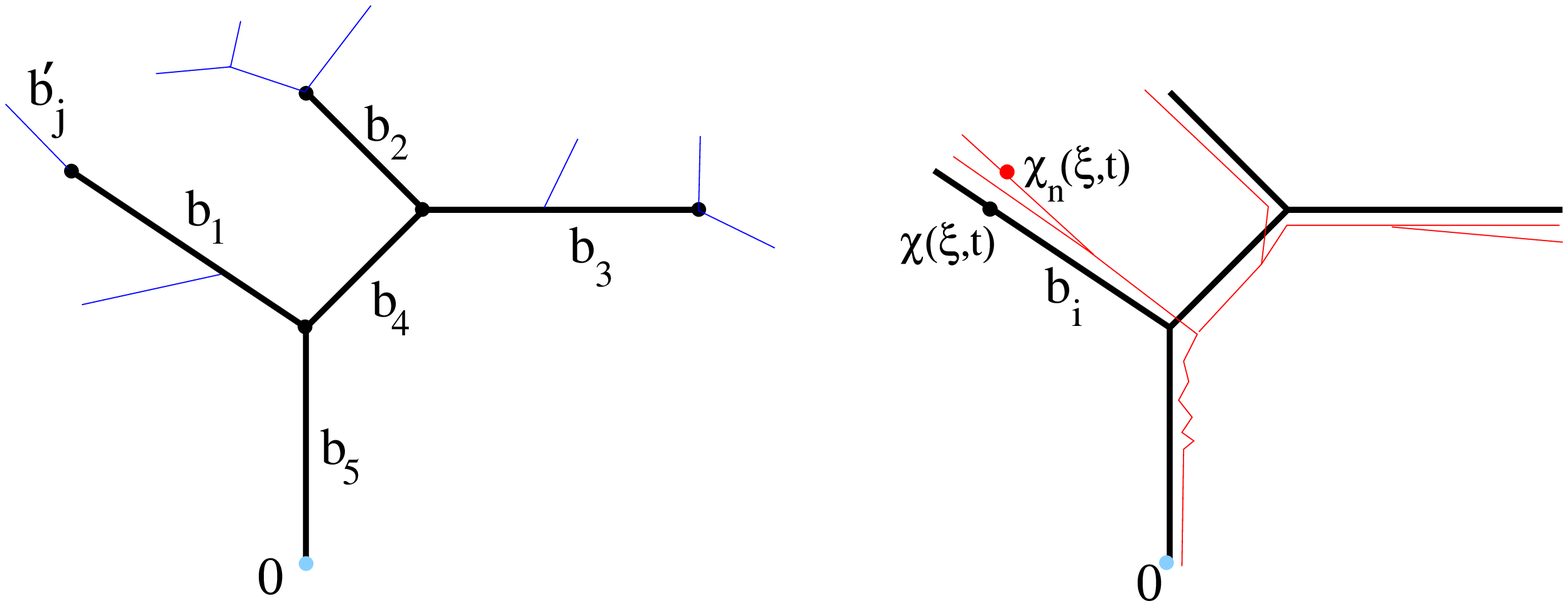}}
\caption{\small Left: Two finite trees, showing three maximal $\ve$-good paths (thick lines) and  8 
maximal $\ve'$-good paths (thin lines), for $0<\ve'<\ve$.  Right: proving the lower semicontinuity of the weighted irrigation cost. 
Given a sequence of irrigation plans $\chi_n\to\chi$,
one can compare
the  cost of $\chi$ restricted to each branch $b_i$ with multiplicity $m(\xi,t)\geq \ve$
with the corresponding costs for the approximating irrigation plans $\chi_n$. }
\label{f:ir86}
\end{figure}

\begin{remark}\label{wemon}  {\rm In the next section we will prove that
 \bel{mp}
\ve'\,<\,\ve\qquad\implies\qquad W^\ve(\xi,t)~\leq~W^{\ve'}(\xi,t).\eeq
Hence the approximations $W^\ve$ depend monotonically on $\ve$.
As a consequence, we can equivalently write
\bel{Wlim}
W(\xi,t)~=~\lim_{\ve \to 0+} W^\ve(\xi,t).\eeq
One should be aware  that this limit may well be $+\infty$.   
}
\end{remark}

\begin{remark} {\rm The assumption {\bf (A2)}, introduced below (\ref{mineq}), guarantees that
the approximation is meaningful.   To see what goes wrong when  {\bf (A2)} fails, 
consider the irrigation plan $\chi:[0,M]\times\R_+\mapsto\R^2$ defined as
$$\chi(\xi,t)~=~\left\{\bega{rl} (t\xi, t) \qquad &\hbox{if}~~ t\in [0,1],\\[3mm]
(\xi,1)\qquad &\hbox{if}~~ t\geq 1.\enda\right.$$
In this case the multiplicity is $m(\xi,t)=0$ for all $\xi\in [0,M]$ and  $t>0$.
Hence $W^\ve(\xi,t)\equiv 0$ for all $\ve>0$.
}\end{remark}
\v
Having constructed a family of weights $W(\xi,t)$, we can now define the corresponding irrigation cost.  Instead of the function $\psi(s)=s^\alpha$
with  $0<\alpha\leq 1$, 
one can here consider more general cost functions $\psi:\R_+\mapsto\R_+$,
satisfying the same assumptions imposed on $f$ at (\ref{fp1}).
As usual, an upper dot will denote a derivative w.r.t.~time.

\begin{definition} Let $f,\psi:\R_+\mapsto\R_+$ be continuous functions, both satisfying
all the assumptions in  {\bf (A1)}.
Let $\chi $ be an irrigation plan satisfying {\bf (A2)} and let $W=W(\xi,t)$ be the 
corresponding weight function, as in (\ref{W}).
If each path is parameterized by arc-length, the {\bf  weighted cost} is then defined as
\bel{EW} \E^{W,\psi}(\chi)~\doteq~\int_0^M\int_0^{\tau(\xi)} 
{\psi(W(\xi,t)) \over m(\xi,t)} \,dt\, d\xi\,.\eeq
More generally, for an arbitrary parameterization of the paths $\chi(\xi,\cdot)$, the 
weighted cost is
\bel{EWP} \E^{W,\psi}(\chi)~\doteq~\int_0^M\int_0^{\tau(\xi)} 
{\psi\bigl(W(\xi,t)\bigr) \over m(\xi,t)} \, |\dot\chi(\xi,t)|\,dt\, d\xi\,.\eeq
\end{definition} 
\v
\begin{remark} {\rm In the special case where $f \equiv 0$, the weight function coincides with
the multiplicity:
$W(\xi,t)=m(\xi,t)$.  Taking $\psi(s)=s^\alpha$ for some $0<\alpha\leq 1$, by (\ref{mineq}), this implies
$\E^{W,\psi}(\chi)\geq\E^\alpha(\chi)$.  Equality holds whenever $\chi$ has the single path property and hence $m(\xi,t)=|\chi(\xi,t)|$.
}
\end{remark}

In order to compute an approximate value of the weighted cost, fix any $\ve>0$
and
let $\Hat \gamma_1,\ldots \Hat\gamma_\nu$ be the maximal $\ve$-good paths.
Consider the elementary paths $\gamma_i$ constructed by the path splitting algorithm 
{\bf (PSA)} at (\ref{eph}),
and let $W_i^\ve:[a_i, b_i]\mapsto [\ve, +\infty[$ be the corresponding approximate weights constructed at (\ref{W1})--(\ref{posti}).  We claim that
\bel{cv1}\E^{W^\ve,\psi}(\chi)~\doteq~\int_0^M \int_0^{\tau(\xi)}
{\psi(W^\ve(\xi,t))\over m(\xi,t)}|\dot{\chi}(\xi,t)|\, dt\, d\xi~=~\sum_{i=1}^N
\int_{a_i}^{b_i} \psi(W^\ve_i(s))\, ds .\eeq

Indeed, recalling (\ref{tep}), denote by $\Omega_\ve\subseteq [0, M]$ the set of particles such that $\tau_{\ve}(\xi) > 0$. By the definition of approximate 
weights $W^\ve$ at (\ref{We1}), 
it follows 
\bel{2a}\int_0^M \int_0^{\tau(\xi)}
{\psi(W^\ve(\xi,t))\over m(\xi,t)}|\dot{\chi}(\xi,t)|\, dt\, d\xi~=~\int_{\Omega_\ve} \int_0^{\tau_{\ve}(\xi)}{\psi(W^\ve(\xi,t))\over m(\xi,t)}|\dot{\chi}(\xi,t)|\, dt\, d\xi.  \eeq
For each $\xi\in \Omega_\ve$, define
$$s_\ve(\xi)~\doteq ~\int_0^{\tau_{\ve}(\xi)} |\dot{\chi}(\xi,t)|\, dt.$$
To fix the ideas, assume that  $\chi(\xi,\cdot)\Big|_{[0, \tau_\ve(\xi)]}\simeq \Hat \gamma_i\Big|_{[0, s_\ve(\xi)]}$ 
for some maximal $\ve$-good path $\Hat \gamma_i$.
Recalling (\ref{We1}), 
 by a standard change of variable formula we obtain
\bel{4a}\int_0^{\tau_{\ve}(\xi)}{\psi(W^\ve(\xi,t))\over m(\xi,t)}|\dot{\chi}(\xi,t)|\, dt~=~\int_0^{s_\ve(\xi)}{\psi(\Hat W_i(s)) \over \Hat m_i(s) }\, ds .\eeq
For each $s>0$ 
consider the set
\bel{5a}\Omega_{i,k}(s)~\doteq~\Big\{  \xi\in [0,\, M];~
\chi(\xi,\cdot)\Big|_{[0,t]}\simeq \Hat \gamma_i\Big|_{[0, s ]} \hbox{ for some }
t>0,~t_{i,k-1}< s\leq t_{i,k}\Big\} . \eeq
Splitting the integral in (\ref{4a}) over the disjoint
intervals $]t_{i, k-1}, t_{i,k}]$ considered at (\ref{tij}),   one obtains 
\bel{6a} \int_0^{s_\ve(\xi)}{\psi(\Hat W_i(s)) \over \Hat m_i(s) }\, ds
~=~\sum_k \int_{t_{i,k-1}}^{t_{i,k}}\Big[ {\psi(\Hat W_i(s))\over \Hat 
m_i(s)}\mathbf{I}_{\Omega_{i,k}(s)}(\xi)\Big]\, ds ,\eeq
where $\mathbf{I}_{\Omega_{i,k}(s)}(\cdot) $ is the indicator function of set 
$\Omega_{i,k}(s)$.
Observing that 
$$\int_{\Omega_\ve}\Big[ {\psi(\Hat W_i(s))\over \Hat 
m_i(s)}\mathbf{I}_{\Omega_{i,k}(s)}(\xi)\Big]\, d\xi~=~{ \psi(\Hat W_i(s)) }, $$
we eventually obtain (\ref{cv1}).

\v
The next lemma shows that the family of approximating weight functions $W^\ve$
is monotonically increasing as $\ve\downarrow 0$.
\begin{lemma} Let $\chi$ be an irrigation plan and let 
the approximate weights $W^\ve$ 
	be defined  as in  (\ref{W1})-(\ref{posti}). Then  for any $0<\ve'<\ve$ 
	and $0\leq s< t$ one has
	\bel{W12}
	W^{\ve}(\xi,t)~\leq~W^{\ve'}(\xi,t).\eeq
\end{lemma}
\v
 {\bf Proof.} 
To prove (\ref{W12}), let $\ve'<\ve$ and let $\tau_{\ve'}(\xi)\geq\tau_\ve(\xi)$  be the
corresponding stopping times in (\ref{tep}).  
By construction, it trivially follows
\bel{sos1}W^{\ve}(\xi,t)~=~W^{\ve'}(\xi,t)~=~0\qquad\forall t\geq \tau_{\ve'}(\xi),\eeq
\bel{sos2}W^\ve(\xi,t)~=~0~\leq~W^{\ve'}(\xi,t)\qquad\forall t\in\, ] \tau_\ve(\xi),\,\tau_{\ve'}(\xi)].\eeq

To prove the inequality in (\ref{W12}) for $t\leq\tau_\ve(\xi)$,
let $\Hat \gamma_1',\ldots, \Hat \gamma'_{\nu'}$ be maximal $\ve'$-good
paths, and let $\gamma'_1,\ldots, \gamma'_{N'}$ be the corresponding
elementary paths, generated by the algorithm {\bf (PSA)}.    
By definition, the weights $W^{\ve'}$ are obtained by induction, performing the 
steps (i)-(ii) at (\ref{W1})--(\ref{posti}) for the elementary paths $\gamma_i'$.

Consider the functions
$$f^\ve_i(w,s)~=~\left\{ \bega{cl} f(w)\quad &\hbox{if}\quad m_i(s)\geq\ve,\\[3mm]
0\quad &\hbox{if}\quad m_i(s)<\ve.\enda\right.$$
Performing the same inductive construction, but with $f$ replaced by $f_i^\ve$
on each elementary path $\gamma_i'$, $i=1,\ldots, N'$, we now recover
exactly the weights $W^\ve$.  A comparison argument now yields 
(\ref{W12}), for all $\xi,t$.
\endproof

As a consequence, we have
\begin{corollary}\label{l:We1} Let $\chi$ be an irrigation plan which satisfies the 
assumption {\bf (A2)}.
Then the weighted irrigation cost in (\ref{EW}) is computed by
\bel{WC4}
\E^{W,\psi}(\chi)~=~~\lim_{\ve\to 0+}\, \int_0^M \int_0^{\tau_\ve(\xi)}{\psi(W^\ve(\xi,t))\over m(\xi,t)}|\dot{\chi}(\xi,t)|
\, dtd\xi
\,.\eeq
\end{corollary}

\section{Lower semicontinuity}
\label{s:6}
\setcounter{equation}{0}

The goal of  this section is to establish the lower semicontinuity of the weighted cost functional 
$\E^{W,\psi}(\chi)$ w.r.t.~pointwise convergence of the irrigation plans.

More precisely, consider a sequence of irrigation plans $\chi_n:[0,M]\times\R_+\mapsto\R^d$.
We say that $\chi_n\to \chi$ {\em pointwise}  
if, for a.e.~$\xi\in [0,M]$, as $n\to \infty$
one has the convergence 
\bel{pconv}
\chi_n(\xi,t)~\to~ \chi(\xi,t)\qquad\hbox{uniformly for $t$ in compact intervals.}\eeq
In terms of the distance (\ref{dgg}), this means
$$\lim_{n\to\infty}~d\bigl(\chi_n(\xi,\cdot),\,\chi(\xi,\cdot)\bigr)~=~0
\qquad\hbox{for a.e. }~~\xi\in [0,M].$$

\begin{theorem}\label{t:51} Consider a sequence $(\chi_n)_{n\geq 1}$ of irrigation plans,
all satisfying the assumption {\bf (A2)}, 
pointwise converging to an irrigation plan $\chi$. 
Assume that the functions $f,\psi$ both satisfy the conditions in {\bf (A1)}.
Then the corresponding
weighted costs satisfy
\bel{lsc}
 \E^{W,\psi}(\chi)~\leq~\liminf_{n\to \infty} \E^{W,\psi}(\chi_n).\eeq
\end{theorem}

\begin{figure}[ht]
\c{\includegraphics[scale=0.6]{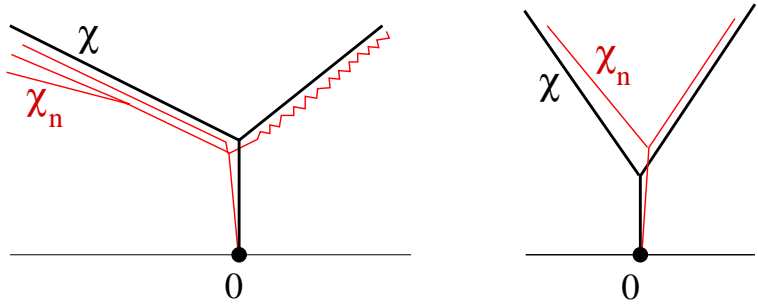}}
\caption{\small Left: two cases where the inequality (\ref{lsc}) can be strict. (i)
Paths in $\chi_n$ may remain separate, while in $\chi$ they all join together.
(ii) Paths in $\chi_n$ may converge to a  path in $\chi$ with strictly smaller length.
Right: two cases  where the weighted irrigation costs satisfy 
$\E^{W,\psi}(\chi_n)<\E^{W,\psi}(\chi)$.
(i) The paths in $\chi_n$ can be slightly shorter than those in $\chi$.
(ii) Paths in $\chi_n$ may remain joined together for a slightly longer time than those in $\chi$.  However, these differences vanish asymptotically, as $n\to\infty$. }
\label{f:ir90}
\end{figure}

Toward a proof, some preliminary results will be needed.

\begin{lemma}\label{l:pathl}
Let $\gamma:[a,b]\mapsto\R^d$ be a Lipschitz path, and let $\ve>0$.
Then there exists $\delta>0$ such that, for any Lipschitz
path $\gamma^\dagger:[a,b]\mapsto\R^d$
which satisfies
$$|\gamma^\dagger(s)-\gamma(s)|~\leq~\delta
\qquad\qquad \forall s\in [a+\delta,\,b-\delta],$$
the length of $\gamma^\dagger$ is bounded below by
\bel{ggd}
\int_{a+\delta}^{b-\delta} |\dot \gamma^\dagger(s)|\, ds~\geq~
(1-\ve)\int_a^b |\dot \gamma(s)|\, ds.\eeq
\end{lemma}

{\bf Proof.} This is an immediate consequence of the lower semicontinuity of the 
path length.\endproof 
\v
In the forthcoming analysis, it will be convenient to use
a distance between two paths which is independent of their parameterization. 
For this purpose, following  \cite{BR} we introduce

\begin{definition}  \label{d:pfm}{\bf (Parameterization-free distance among paths).}
	Given two continuous paths $\varphi_i: [0,\, S_i]\mapsto \R^d$, $i = 1,2$, the distance $\delta(\varphi_1,\varphi_2)$ is defined as
	\bel{2e}\delta(\varphi_1,\varphi_2)~\doteq~\inf_{\eta_1,\eta_2}\max_{t\in [0,1]}
	\Big|\varphi_1(\eta_1(t)  - \varphi_2(\eta_2(t)) \Big|\, ,\eeq
	where the infimum is taken over all couples of continuous, nondecreasing, surjective maps $\eta_i:[0,1]\mapsto [0,\, S_i].$
\end{definition}

\begin{figure}[ht]
\c{\includegraphics[scale=0.5]{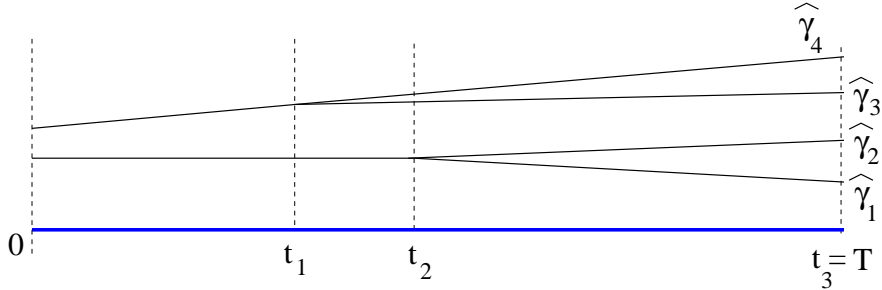}}
\caption{\small  The two configurations compared in Lemma~\ref{l:mpc}.
For every $t\in [0,T]$, the sum of the weight functions $W_i(t)$ along a family of maximal paths
$\Hat \gamma_i$ 
is compared with a single weight $W(t)$, satisfying the ODE (\ref{WTT}).}
\label{f:ir96}
\end{figure}

As shown in \cite{BR},  one has
\begin{itemize}
		\item[i)] $\delta(\varphi_1,\varphi_2)~=~\delta(\varphi_2,\varphi_1)~\geq~0\, ,$
		\item[ii)] $\delta(\varphi_1,\varphi_2)~=~0 \hbox{ if and only if } \varphi_1~\simeq~\varphi_2,$ in the sense of Definition~\ref{d:eq1}\, ,
		\item[iii)] $\delta(\varphi_1,\varphi_3)~\leq~\delta(\varphi_1,\varphi_2) + \delta(\varphi_2,\varphi_3)\, .$ 
\end{itemize}
The proof of the following  lemma is elementary, but the conclusion turns out to be crucial in the proof of lower semicontinuity of the irrigation cost.
\begin{lemma}\label{l:doc1}
	Let $\gamma_i:[0,\, \ell_i]\mapsto \R^d$, $i=1,2$, be two 
	paths parametrized by arc-length.
	Assume that they bifurcate at some time $0\leq \tau<\min\,\{\ell_1,\ell_2\}$, i.e.
$$\tau~=~ \sup\Big\{ t\geq 0\,;~~\gamma_1(s)=\gamma_2(s)\quad\forall s\in [0,t]\Big\}.$$
Then for any $h>0$, there exists $\sigma>0$ such that
	\bel{4e}\delta\Big(\gamma_1\Big|_{[0,s]},\gamma_2\Big|_{[0,t]}\Big)~
	\geq~\sigma,\qquad \forall s\in [\tau+h,\, \ell_1],\,  t\in[0,\, \ell_2] \, . \eeq
\end{lemma}
{\bf Proof.} The map $(s,t)\mapsto 
\delta\bigl(\gamma_1|_{[0,s]},
\,\gamma_2|_{[0,t]} \bigr)$ is continuous and strictly positive on the compact 
domain 
$[\tau+h,\, \ell_1 ]\times [0,\, \ell_2 ]$.  Hence it has a strictly positive minimum.
\endproof
\v
In Lemma~\ref{l:simple} we compared the weight $w(s)$ along a single path 
$\gamma$ with a sum of weights
$\sum_i w_i(s)$ along a family of distinct paths $\gamma_i$.  
The next lemma extends this result to 
a more general configuration where the paths $\gamma_i$  are not necessarily disjoint, as shown in 
Fig.~\ref{f:ir96}.

More precisely, consider
an irrigation plan $\chi$ containing finitely many  maximal paths $\Hat\gamma_j:[0,T]\mapsto \R^d$, $ j = 1,\ldots, \nu$, all parameterized by arc-length and all with the same length $T$.
Let $\Hat m_j:[0,T]\mapsto \R^d$ be the (non-increasing)
multiplicity function along $\Hat\gamma_j$, defined as in (\ref{mi}), and consider weights
\bel{TVJ}
\Hat W_j (T) ~\geq ~\Hat m_j(T)~>~0,\eeq
arbitrarily assigned at the terminal point of each maximal path.
In turn, these data 
determine the weight functions along all paths.   Namely, 
let $\gamma_i:[a_i, b_i]\mapsto \R^d, 1\leq i\leq N$ be the corresponding elementary paths, constructed by the Path Splitting Algorithm {\bf (PSA)}.
By backward induction we can now construct the weights $W_i$ along each elementary path,
in a similar way as in (\ref{W1})--(\ref{posti}).
 \begi
\item  For every index $i$ such that $b_i=T$, the weight $W_i :[a_i, b_i]\mapsto \R_+$ along 
the elementary path $\gamma_i$ is computed by solving
 \bel{WIT}
w(t)~=~\int_{t}^{b_i} f(w(s))\, ds + m_i(t) + \bigl[ \Hat W_{j(i)}(T) -  m_i (T) \bigr]\,,
\qquad\qquad t\in \,]a_i, T].
\eeq
Here $\Hat \gamma_{j(i)}$ is the unique maximal path that contains $\gamma_i$ as its restriction 
to $[a_i, b_i] = [a_i,T]$. 
\item If $b_i<T$,
the weight $W_i:[a_i, b_i]\mapsto\R_+$ along the 
elementary path $\gamma_i$  is then defined to be the solution of
 \bel{Wij}
w(t)~=~\int_{t}^{b_i} f(w(s))\, ds + m_i(t) + \left[\sum_{k\in \O(i)} W_k(a_k+) -\sum_{k\in \O(i)} m_k(a_k+)\right],\qquad t\in \,]a_i, b_i].
\eeq
As in (\ref{posti}), here the summations range over all elementary paths $\gamma_k$ 
that originate from the tip of $\gamma_i$.
\endi

\begin{lemma}\label{l:mpc}  Let the weights $W_i:[a_i, b_i]\mapsto\R_+$ be constructed as above.
Given any constant $\Hat W$ such that
 \bel{0050} 0~<~ \Hat W~\leq~\sum_{j=1}^\nu  \Hat W_j(T)\,,   \eeq 
let $W:[0,T]\mapsto\R$ be the solution to  the backward Cauchy problem
\bel{WTT}\dot W(t)~=~ - f(W(t)), \qquad\qquad W(T)\,=\, \Hat W.\eeq
Then for all $t\in \,]0,T]$ one has
\bel{0054} W(t)~\leq~\sum_{i\in I(t)} W_i(t)\,, \eeq
where $I(t)$ denotes the set of indices $i\in \{1,\ldots,N\}$ such that $a_i<t\leq b_i$. 
As a consequence, 
\bel{1004} \int_0^T \psi(W(t))\, dt~\leq~\sum_{i=1}^N \int_{a_i}^{b_i} \psi(W_i(t))\, dt\,.   \eeq
\end{lemma}

{\bf Proof.} Let $0<t_1<\cdots<t_q=T$ be the times where two or more maximal paths bifurcate.
The proof will achieved by backward induction on $p=1,2,\ldots,q$.
\v
{\bf 1.}
For $t\in ]t_{q-1}, t_q] = \, ]t_{q-1}, T]$, the above definition implies $I(t)= I(T)$. 
By (\ref{0050})  it follows
\bel{1054} \sum_{i\in I(T)} W_i(T)~\geq~\Hat W\, .   \eeq
 For each $i\in I(T)$, $t\in \,]t_{q-1}, T]$, by (\ref{WIT}) it follows
\bel{0056} W_i(t)~=~\int_t^{t_q} f(W_i(s))\, ds + W_i(T) + \left[m_i(t) - m_i(T)\right]   \eeq
On the other hand,
\bel{1002} W(t)~=~\int_t^{t_q} f( W(s))\, ds +\Hat W\, .  \eeq
Because of (\ref{1054}) we can apply Lemma~\ref{l:simple} and conclude that
\bel{0055} \sum_{i\in \I(T)} W_i(t)~\geq~ W(t),\qquad \forall t\in \,]t_{q-1}, T ].   \eeq
\v
{\bf 2.}
Next, assume that that the inequality (\ref{0054})  has been proved for all 
$t\in \,]t_p, T]$,
for some $1\leq p <q$. We claim that it also holds for $t\in \,]t_{p-1}, t_p]$. 

Indeed, since the solution $W$ of (\ref{WTT}) is continuous while all weights $W_i$ 
are non-increasing, the inductive assumption yields
\bel{0072}  \sum_{i\in \I(t_{p})} W_i(t_{p})~\geq~ W(t_p)\, .  \eeq
For each $i\in I(t_{p})$, $t\in \,]t_{p-1}, t_{p}]$, we then have
\bel{0071}\bega{c}\ds W_i(t)~=~\int^{t_{p}}_t f(W_i(s))\, ds + W_i(t_{p}) + \left[ m_i(t) - m_i(t_{p})  \right]\,,\\[3mm]
\ds W(t)~=~\int_t^{t_{p}} f( W(s))\, ds +  W(t_{p})  \enda\eeq
Because of  (\ref{0072}) we can again apply Lemma~\ref{l:simple} and conclude
\bel{0073} \sum_{i\in \I(t_{p})} W_i(t)~\geq~ W(t),\qquad \forall t \in (t_{p-1},\, t_{p}].   \eeq
By induction on $p$, this yields a proof of (\ref{0054}). 
\v
{\bf 3.}
Since $\psi$ satisfies the assumption {\bf (A2)}, from (\ref{0054}) it follows
\bel{1003}\bega{l}\ds \sum_{i=1}^N\int_{a_i}^{b_i} \psi(W_i(t))\, dt~=~\sum_{p=1}^q\sum_{i\in I(t_p)} \int_{t_{p-1}}^{t_p} \psi(W_i(t))\, dt \\[3mm]
\ds \quad\geq~\sum_{p=1}^q \int_{t_{p-1}}^{t_p} \psi( W(t))\, dt~=~\int_0^T \psi(W(t))\, dt  \enda\eeq
Hence  (\ref{1004}) holds.
\endproof

\begin{remark}\label{r:wcomp}{\rm
In Lemma~\ref{l:mpc} we assumed that all maximal paths $\Hat\gamma_j$ had the same length $T$.
The same conclusions (\ref{0054})-(\ref{1004}) remain valid if each maximal path $\Hat\gamma_j:
[0, T_j]\mapsto\R^d$ is defined on an interval of length $T_j\geq T$, replacing (\ref{0050})
with
\bel{51}
\Hat W~\leq~\sum_{j=1}^\nu \Hat W_j(T_j).\eeq  To prove this, it suffices to 
consider the restriction of each $\Hat\gamma_j$ to the sub-interval $[0,T]$, and observe that
(\ref{51}) implies (\ref{0050}), because the weight functions are non-increasing, 
}\end{remark}

\v
After these preliminaries, we are now ready to give a proof of the main result of this section, in several steps.

{\bf Proof of Theorem~\ref{t:51}.} 

{\bf 1.} Without loss of generality, we can assume that all paths $\chi_n(\xi,\cdot)$ are parameterized by arc-length. 
As a consequence,  for each $\xi\in [0,\, M]$, the limit paths $\chi(\xi,\cdot)$ 
will be 1-Lipschitz, but not necessarily parameterized by arc-length.   

Fix any $\ve_0>0$.
Let $\tau_{\ve_0}$ be the corresponding stopping time as in (\ref{tep}),
and define the  truncated irrigation plan  
\bel{trc}\chi^{\ve_0}(\xi,t)~\doteq~\left\{\bega{cl}\chi(\xi,t)\qquad &\hbox{if}\qquad t\leq \tau_{\ve_0}(\xi),\\[4mm]
\chi(\xi,\tau_{\ve_0}(\xi))\qquad &\hbox{if}\qquad t>\tau_{\ve_0}(\xi).\enda\right. \eeq
Using Corollary~\ref{l:We1}, the theorem will be proved by showing that
\bel{ntry1} 
\E^{W,\psi}(\chi^{\ve_0})~\leq~\liminf_{n\to\infty} \E^{W,\psi}(\chi_n) .\eeq


\v
{\bf 2.} 
For each $\xi\in [0,\, M]$, in order to re-parameterize the limit 
path $\chi(\xi,\cdot)$ in terms of arc-length, let 
\bel{p0p}s(\xi,t)~\doteq~ \int_0^t|\dot{\chi}(\xi,r)|\, dr.\eeq
A left-continuous 
inverse of $s(\xi,\cdot)$, taking values in  $\R_+\cup\{+\infty\}$, can be defined as 
\bel{1k}\eta(\xi,s)~\doteq~\inf\Big\{t~\geq~0;~s(\xi,t)~=~s\Big\}\, .\eeq
The map 
\bel{lcrp}
s~\mapsto~\chi(\xi, \eta(\xi,s))\eeq
now provides the arc-length parameterization of $\chi(\xi,\cdot)$.
We observe that, for each $s$, the map  $\xi\mapsto \eta(\xi,s)$ 
 is measurable. Moreover, since $|\dot{\chi}(\xi,t)|\leq 1$, one has 
\bel{etaL}\eta(\xi, s_2)-\eta(\xi,s_1)~\geq ~s_2-s_1\qquad\forall 0\leq s_1<s_2\,.\eeq
\v
{\bf 3.} 
Next, let $\Hat \gamma_1,\ldots,\Hat \gamma_{\nu}$ be the maximal $\ve_0$-good 
paths for the irrigation plan $\chi$. As before, we assume that each
 $\Hat\gamma_j:[0,\Hat \ell_j]\mapsto\R^d$ is parameterized by arc-length.
For $s\in\, ]0, \Hat{\ell}_j]$, let 
$$\Hat\Omega_j(s)~\doteq~\left\{\xi\in [0,M];\quad 
\chi(\xi,\cdot)\Big|_{[0,t]} \simeq \Hat\gamma_j\Big|_{[0,s]} ~\hbox{ for some } t> 0  \right\}$$
be the set of particles whose trajectory follows the path $\Hat\gamma_j$, at least up to the point $\Hat \gamma_j(s)$.

Implementing  the algorithm {\bf (PSA)} 
described 
at (\ref{tij})--(\ref{eph}), these maximal paths can be 
split into finitely many elementary paths
$\gamma_1,\ldots, \gamma_N$.  By construction, each
 $\gamma_i:[a_i, b_i]\mapsto \R^d$, $1\leq i\leq N$ is the restriction of some $\Hat \gamma_j$   to a subinterval $[a_i,b_i]$.
 We then define 
\bel{+4}\Omega_i(s)~\doteq~\Hat\Omega_j(s)\qquad \forall s\in\, [a_i,b_i ]. \eeq
The multiplicity function $m_i: [a_i, b_i]\mapsto \R_+$ along the 
elementary path $\gamma_i$ is then computed by
\bel{+5} m_i(s) ~=~ \meas(\Omega_i(s))\, ,\qquad s\in\, [a_i, b_i]\, . \eeq
According to (\ref{cv1}), the approximate weighted irrigation cost is computed by
a sum over all elementary paths:
\bel{EW0}
E^{W,\psi}(\chi^{\ve_0})~=~\sum_{i=1}^N \int_{a_i}^{b_i} \psi(W_i^{\ve_0} (s))\, ds,\eeq
where the weights $W^{\ve_0}_i$ are determined as in (\ref{W1})--(\ref{posti}).
\v
\begin{figure}[ht]
\c{\includegraphics[scale=0.4]{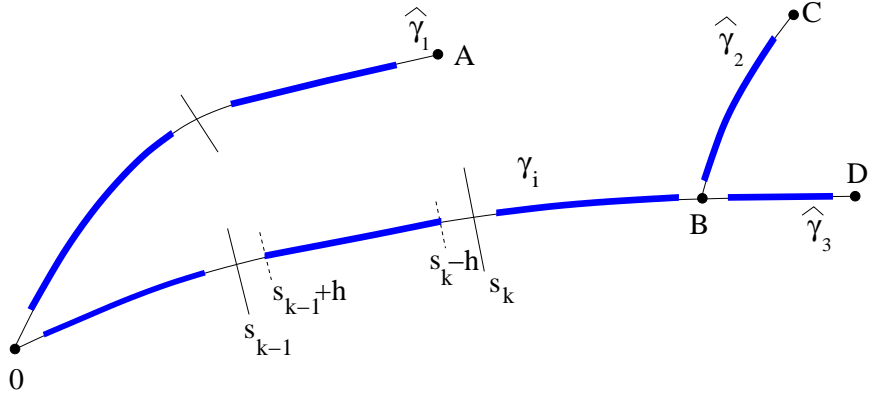}}
\caption{\small  Proving
the lower semicontinuity of the weighted irrigation cost, steps 3-4.
Here $\Hat\gamma_1,\Hat\gamma_2,\Hat\gamma_3$ are maximal $\ve_0$-good 
paths of $\chi$, while $\gamma_i:[a_i, b_i]\mapsto\R^d$, with endpoints $0,B$, 
is an elementary  path produced by the algorithm {\bf (PSA)}.
The  path $\gamma_i$ 
is further partitioned, taking subintervals $[s_{k-1}, s_k]$ of length $\leq\delta$.
We then approximate
the multiplicity $m_i(s)$ with a piecewise  constant 
function $\Tilde m_i$, as in (\ref{TMI})  and replace $f$ with $f^h$ as in (\ref{TMI1}).
By choosing the constants $\delta,\delta_0,h>0$ sufficiently small, the  new weight
$\Tilde W_i$ determined by (\ref{W1T}) can be kept arbitrarily close to the original weight $W_i^{\ve_0}$.  
 }
\label{f:ir95}
\end{figure}

{\bf 4.} 
We claim that it is possible to replace the multiplicity functions 
$m_i$ by  strictly smaller piecewise constant functions $\Tilde m_i$, 
producing a very small change in the weights $W^{\ve_0}_i$.  
More precisely, for each $i\in\{1,\ldots,N\}$,  choose $\delta>0$ and insert 
the times (see Fig.~\ref{f:ir95})
\bel{+0}a_i\,=\,s_0\,<\,s_1\,<\,\cdots\,< \,s_{n(i)} \,=\,b_i\,,\eeq
so that $s_k - s_{k-1}\leq \delta$ for every $k=1,\ldots, n(i)$.
For a given $\delta_0>0$, with $\delta_0<\ve_0$, we  then define the piecewise constant function
\bel{TMI}
\Tilde m_i(t)~=~m_i(s_k)-\delta_0\qquad\forall t\in\,]s_{k-1}, s_k].\eeq
Since $s\mapsto m_i(s)\in [\ve_0, +\infty[\,$ is non-increasing, we clearly have 
$0<\Tilde m_i(t)< m(t)$ for all $t\in \,]a_i, b_i]$.  
Next, given another  constant $h>0$, with $h<\!<\delta$, we define
\bel{TMI1} f^{h}(t,\omega)~\doteq~\left\{\bega{cl} f(\omega)&\qquad 
\hbox{if}~~~t\in [s_{k-1} + h,~ s_k - h],\quad1\leq k\leq n(i),\\[3mm]
0&\qquad \hbox{otherwise.}  \enda\right. \eeq

We claim that, for any $\ve>0$, one can choose the above  constants
$\delta, \delta_0, h>0$ small enough so that, replacing the multiplicities $m_i$ with $\Tilde m_i$, and replacing $f$ with $f^{h}$, 
the corresponding weights $\Tilde W_i$ satisfy
\bel{WTW}
 \|W_i^{\ve_0}-\Tilde W_i\|_{\L^1([a_i,b_i])}~<~\ve,\qquad 
\qquad
|W_i^{\ve_0}(a_i+)-\Tilde  W_i(a_i+)|\,<\,\ve \,.\eeq

Indeed, recalling (\ref{i1n}) consider first the case $i\in \I_1$, so that $\gamma_i$ is one of the outer-most branches. 
Then the weight 
$W_i^{\ve_0}$ is obtained by solving (\ref{W1}), while  $\Tilde W_i$ provides a solution to 
\bel{W1T}
w(t)~=~\int_t^{b_i}f^h(s,w(s))\, ds + \Tilde m_i(t),\qquad\qquad t\in \,]a_i, b_i].
\eeq
By choosing $\delta, \delta_0, h >0$ sufficiently small, we can make
the differences 
$$\|m_i - \Tilde m_i\|_{\L^1([a_i,b_i])}\,,\qquad |m_i(a_i+) - \Tilde m_i(a_i+)|\,,\qquad \meas\left([a_i,\, b_i]\setminus \bigcup_{k=1}^{n(i)} [s_{k-1}+h,\,s_k - h] \right)$$
as small as we like.
The estimate (\ref{WTW}) thus follows from part (iii) of Lemma~\ref{l:mvode}.

The case $i\in \I_p$ for $p>1$ is proved in the same way, by induction on $p$.
\v
In view of (\ref{ntry1}) and (\ref{EW0}), to prove the theorem 
it thus suffices to show that, for any given
$\delta,\delta_0,h>0$, the corresponding weights $\Tilde W_i$ satisfy
\bel{lsc2}
\sum_{i=1}^N \int_{a_i}^{b_i} \psi(\Tilde W_i (s))\, ds~\leq~\liminf_{n\to \infty} \E^{W,\psi}(\chi_n).\eeq
\v
{\bf 5.} Consider again the arrival time  $\xi\mapsto \tau(\xi)$ introduced 
	in Definition~\ref{d:27}.  For any $\ve>0$, by Corollary~\ref{c:210} there is
a  compact set $\Omega_{\ve}\subseteq [0,M]$, with 
\bel{mok}\meas\,\bigl([0,M]\setminus \Omega_{\ve} \bigr)~<~\ve\,,\eeq 
on which that map $\tau(\cdot)$ is continuous.   Hence
	\bel{+1}\max_{\xi \in \Omega_{\ve}}~ \tau(\xi)~\leq~\kappa \eeq
	for some constant $\kappa$. By (\ref{pconv}) and Egoroff's theorem, by 
slightly shrinking the compact set $\Omega_{\ve}$, we can assume that 
(\ref{mok}) still holds, together with
	\bel{+2} \lim_{n\to \infty} \sup_{\xi\in \Omega_{\ve}}\|\chi(\xi,\cdot) - \chi_n(\xi,\cdot) \|_{\mathbf{L}^\infty([0,\kappa])}~=~0\,. \eeq
In addition, calling $\tau^n(\xi)$ the smallest time $\tau$ such that $\chi_n(\xi,\cdot)$ is constant for $t\geq \tau$, 
 by further shrinking  $\Omega_{\ve}$ we can also assume
\bel{+18}\liminf_{n\to \infty} \inf_{\xi\in \Omega_{\ve}} \left[ \tau^n(\xi) - \tau(\xi)  \right]~\geq~0 . \eeq
Indeed, since
\bel{+17}\liminf_{n\to \infty} \tau^n(\xi)~\geq~\tau(\xi),\eeq
it follows that the non-decreasing sequence
$$\Hat\tau^n(\xi)~\doteq~\inf_{k\geq n}\tau^k(\xi)$$
converges to a limit
$$\lim_{n\to\infty}\Hat\tau^n(\xi)~=~\tau^\infty(\xi)~\geq~\tau(\xi)$$
for a.e.~$\xi\in [0,M]$.  Again by Egoroff's theorem we can choose a large subset
$\Omega_\ve\subset [0,M]$
where the pointwise convergence is uniform. This yields (\ref{+18}).

Furthermore, since each $\chi_n$ satisfies the assumption {\bf (A2)}, we can choose
$\ve_n>0$ small enough so that the following holds.  Defining the stopping time
\bel{tepn}
\tau^n_{\ve_n}(\xi)~\doteq~\sup\,\bigl\{ t\geq 0\,;~~m_n(\xi,t)\geq\ve_n
\bigr\},\eeq
by possibly further shrinking the 
set $\Omega_\ve$ in (\ref{mok}) one has
\bel{ten}
\tau^n_{\ve_n}(\xi)~\geq~\tau^n(\xi) - {h\over 2}\qquad\forall \xi\in \Omega_\ve, ~~n\geq 1\,.\eeq
\v
{\bf 6.} Let $\Hat\gamma'_1,\ldots, \Hat\gamma'_{\nu'}$ be the maximal $\ve_n$-good paths in 
$\chi_n$, and let $\gamma'_1,\ldots,\gamma_{N'}'$ be 
the elementary paths constructed by the algorithm  {\bf (PSA)}. 
As in step {\bf 3}, for each $\Hat\gamma_j': [0,\hat\ell_j']\mapsto\R^d$ we define  
\bel{+3}\bega{rl}\Hat\Omega_j'(s)&\ds\doteq~\left\{\xi\in [0,M];~\chi_n(\xi,\cdot)
\Big|_{[0,t]} \simeq \Hat\gamma_j'\Big|_{[0,s]} \hbox{ for some } t>0 \right\}\\[4mm]
&=~\ds\Big\{\xi\in [0, M];~\chi_n(\xi, t) = \Hat\gamma_j'(t), 
~\forall t \in [0, s] \Big\}. \enda\eeq
This is the set of particles whose trajectory follows the maximal path $\Hat \gamma_j'$, at least up to time $s$.
Notice that the last identity holds because $\Hat \gamma_j'$ and  $\chi_n(\xi,\cdot)$ are
both parameterized by arc-length. By construction, each elementary path
$\gamma_i':[a_i', b'_i]\mapsto \R^d$, $1\leq i\leq N'$ is the restriction of some $\Hat\gamma_j'$ to a subinterval $[a'_i, b'_i]$. We then define 
\bel{*0} \Omega'_i(s)~\doteq~\Hat\Omega_j'(s)\qquad \forall s \in\, [a'_i, b'_i]. \eeq 
\v
{\bf 7.}  Now consider a  particle $\xi\in \Omega_i(s_k)\cap \Omega_{\ve}$,
so that the path $t\mapsto \chi(\xi,t)$ reaches the point $\gamma_i(s_k)$ at some time
$t=\eta(\xi, s_k)$.
This implies $\tau(\xi)\geq \eta(\xi, s_k)$.  Hence by (\ref{+18}) we have
$$\tau^n(\xi)~>~\eta(\xi, s_k)-{h\over 2}$$
for all $n$ large enough.   In turn, choosing $\ve_n>0$ sufficiently small, by (\ref{ten}) it follows
$$\tau^n_{\ve_n}(\xi)~\geq~\tau^n(\xi)-{h\over 2}~>~\eta(\xi, s_k)-h
~\geq~\eta(\xi, s_k-h).$$
Otherwise stated, by further 
slightly shrinking the compact set $\Omega_{\ve}$ in (\ref{mok}), for  any $h>0$ 
we can thus achieve the  implication
\bel{xet} \xi~\in~\Omega_i(s_k)\cap \Omega_{\ve}\qquad\implies\qquad \eta(\xi, s_k - h)~<~\tau^n_{\ve_n}(\xi),  \eeq
for all $n$ sufficiently large.
\v
{\bf 8.}
We observe that two particles $\xi,\Tilde\xi$, which have the same trajectory in the irrigation plan $\chi$, may be sent along different paths by the irrigation plan $\chi_n$.
To account for this fact, recalling (\ref{+4})    and (\ref{*0}), for a fixed $n\geq 1$ we define 
\bel{+44} A_i^j(s_k)~\doteq~\left\{ \xi\in \Omega_{\ve}\cap \Omega_i(s_k);~\chi_n(\xi,t)~=
~\Hat\gamma_j'(t),~\forall 0\leq t \leq \eta(\xi, s_k - h)  \right\}.  \eeq
In other words, $A_i^j(s_k)$ is the set of particles $\xi\in \Omega_\ve$ such that:
\begi
\item By the irrigation plan $\chi$ they are moved along the $\ve_0$-good
elementary path $\gamma_i$,
at least up to the point $\gamma_i(s_k)$.
\item By the irrigation plan $\chi_n$ they are moved along the $\ve_n$-good maximal path
$\Hat\gamma_j'$, at least up to point  $ \Hat\gamma_j'\bigl(\eta(\xi,s_k-h)\bigr)$.
\endi
Using Lusin's theorem and by possibly shrinking the compact domain 
$\Omega_{\ve}\subseteq [0,M]$, in addition to (\ref{mok}) we can assume that, 
restricted to each $A_i^j(s_k)$,  the two maps
$$ \eta(\cdot,s_k-h): A_i^j(s_k)\mapsto \R_+,\qquad\qquad
 \eta(\cdot,s_{k-1}+h):A_i^j(s_k)\mapsto \R_+$$
are continuous.    
\v
{\bf 9.} The set of paths 
\bel{gik}\gamma_{i,k}~\doteq~\gamma_i\Big|_{[s_{k-1},\, s_k]}\eeq
comes with an obvious partial ordering.  Namely, we define
\bel{ikprec}
(i,k)~\preceq~(i^\dagger, k^\dagger)\eeq
if the two  elementary paths 
$\gamma_i, \gamma_{i^\dagger}$ for the irrigation plan $\chi$ 
are both contained in some $\ve_0$-good maximal path $\Hat \gamma_j$, and moreover
$s_k\leq s_{k^\dagger}$.

\begin{figure}[ht]
\c{\includegraphics[scale=0.6]{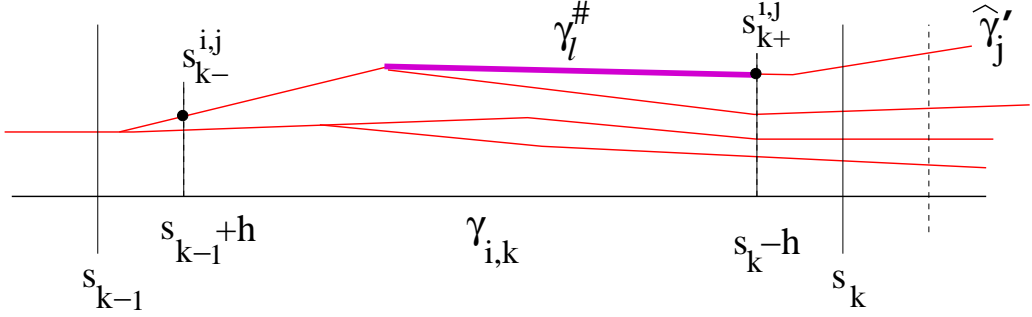}}
\caption{\small  To compare the cost of the irrigation plans $\chi$ and $\chi_n$, to the portion
of the $\ve_0$-good elementary path  $\gamma_i: [s_{k-1}+h, s_k-h]\mapsto\R^d$ we associate
a family of $\ve_n$-good paths $\gamma^\sharp_l$ in $\chi_n$.}
\label{f:ir94}
\end{figure}

As shown in Fig.~\ref{f:ir94}, to each portion $\gamma_{i,k}$
of the elementary path $\gamma_i$ in the irrigation plan $\chi$ we 
shall associate a family $\{\gamma^\sharp_l\}$ of paths  in the irrigation plan $\chi_n$,
and compare the corresponding costs.
For this purpose,  assuming $A_i^j(s_k)\not=\emptyset$, we define
\bel{+11} s_{k+}^{i,j}~\doteq~\inf_{\xi \in A_i^j(s_k)}\eta(\xi,s_k-h),
\qquad\qquad  s_{k-}^{i,j}~\doteq~\inf_{\xi\in A_i^j(s_k)} \eta(\xi, s_{k-1} + h).
\eeq
Notice that, by (\ref{etaL}), one has
\bel{0032} s_{k+}^{i,j} - s_{k-}^{i,j}~\geq~s_k - s_{k-1} - 2h.  \eeq

For each $i,j,k$ such that $A_i^j(s_k)$ is non-empty, we now consider all  the paths
$\gamma^\sharp_{l}:[ a^\sharp_l,\,  b^\sharp_l]\mapsto \R^d$, 
obtained as follows. Consider all the $\ve_n$-good elementary paths $\gamma_p':[a'_p, b'_p]\mapsto \R^d$ of $\chi_n$.   
which are contained in  the maximal path $\Hat\gamma_j'$.   We then take $\gamma_l^\sharp$ 
to be the restriction of 
$\gamma_p'$  
to the subinterval 
\bel{gsd}[ a^\sharp_l,\,  b^\sharp_l]~\doteq~[a'_p, b'_p] \, \cap \,[s_{k-}^{i,j}\,,s_{k+}^{i,j}] .\eeq
Call $\Gamma_{i, k}$ the collection of all  such paths $\gamma^\sharp_l$, 
as $j$ varies
among all the maximal $\ve_n$-good paths of $\chi_n$, with $A_i^j(s_k)\not=\emptyset$.

\v
{\bf 10.}  Let  $m^\sharp_l:[ a^\sharp_l,b^\sharp_l]\mapsto \R_+$ be
the multiplicity of the path $\gamma^\sharp_l$ in the irrigation plan $\chi_n$. 
We claim that,
choosing $\delta_0=\ve$ in (\ref{TMI}), for all $n$ large enough
the piecewise constant multiplicity $\Tilde m_i$ 
defined at (\ref{TMI}) satisfies
\bel{0031} \Tilde m_i(t)~=~\Tilde m_i(s_k - h)~\doteq~m_i(s_k) - \ve~<~\sum_{\gamma_l^\sharp\in \Gamma_{i,k},~b^\sharp_l = s_{k+}^{i,j}} m^\sharp_l(b_l^\sharp)  \eeq 
for all $t\in \,]s_{k-1}, s_k]$.
Indeed, in view of (\ref{xet})-(\ref{+44}),  for each $\xi\in \Omega_{i}(s_k)\cap \Omega_{\ve}$, there is some maximal path $\Hat\gamma_j'$ such that
\bel{4001} \chi_n(\xi,t)~=~\Hat\gamma_j'(t),\qquad \forall t\in [0,\, \eta(\xi, s_k - h)].   \eeq
By (\ref{+11}) we know that $\eta(\xi,s_k - h)\geq s_{k+}^{i,j}$. 
Hence (\ref{4001}) implies $\xi\in \Hat\Omega'_j(s_{k+}^{i,j})$. Therefore
\bel{0035} \Tilde m_i(s_k-h)~=~\meas(\Omega_i(s_k)) - \ve~<~
\meas\Big(\Omega_i(s_k)\cap \Omega_{\ve}\Big)~<~\sum_{\gamma_l^\sharp\in \Gamma_{i,k},~b^\sharp_l = s_{k+}^{i,j}} m^\sharp_l(b_l^\sharp). \eeq

\v{\bf 11.} 
Toward a proof of (\ref{lsc2}),  a key observation is the following.
If two particles $\xi,\Tilde \xi$ are sent  by $\chi$ along two
maximal  paths $\Hat \gamma_i,\Hat\gamma_j$ which bifurcate
at a some time $\tau_{ij}$, then, for all $n$ large enough, the irrigation plan 
$\chi_n$ must send these two particles along distinct paths as well.
In this step we prove a precise estimate in this direction.

Let $\Hat\gamma_1,\ldots, \Hat \gamma_\nu$ be the 
maximal $\ve_0$-good  paths for the
irrigation plan $\chi$.
For any given $h>0$, by Lemma  \ref{l:doc1} 
one can find $\sigma>0$ such that
	\bel{3j}  \delta\left(\Hat\gamma_i\Big|_{[0,s]},\Hat
	\gamma_j\Big|_{[0,t]}\right)~\geq~\sigma\,,\qquad \forall 
	i\neq j\, , ~s\in [\tau_{ij}+h,\,  \hat{\ell}_i],~t\in [0,\, \hat{\ell}_j]\, .  \eeq
	Here $\tau_{ij}$ is the time where the maximal paths $\Hat \gamma_i$ and 
	$\Hat\gamma_j$ bifurcate, as defined at (\ref{tauij}).
	
On the other hand, by (\ref{+2}), for all $n$ sufficiently large one has
\bel{0042} \sup_{\xi \in \Omega_{\ve}} \|\chi(\xi,\cdot) - \chi_n(\xi,\cdot) \|_{\mathbf{L}^\infty([0,\kappa])}~<~{\sigma\over 3}\,,  \eeq
where $\sigma$ is the constant  in (\ref{3j}).

Consider two particles $\xi,\Tilde\xi\in \Omega_\ve$ which are sent by $\chi$ along the two
distinct maximal paths $\Hat \gamma_i,\Hat\gamma_j$.    
More precisely, recalling the definition
(\ref{+3}), assume that for some $h>0$
$$\xi\,\in\, \Hat\Omega_i(\tau_{ij}+h)\cap \Omega_\ve, \qquad\qquad    \Tilde\xi \in \Hat\Omega_j(\tau_{ij}+h)\cap \Omega_\ve\,,
$$ 
Without loss of generality, assume 
$$T~\doteq~ \eta(\xi,\tau_{ij}+h)~\leq ~\eta(\Tilde\xi, \tau_{ij} + h).$$
Recalling the notation used at (\ref{p0p}), we can now find 
$\tau\doteq s(\Tilde\xi, T) \leq \hat{\ell}_j$, such that by (\ref{3j}) and (\ref{0042})
\bel{043}\bega{rl}\sigma&\leq~\ds \delta\Big(\Hat\gamma_i\Big|_{[0,\tau_{ij}+h]}, 
\Hat\gamma_j\Big|_{[0,\tau]}\Big)~=~ \delta\Big(\chi(\xi,\cdot)
\Big|_{[0,T]}, \chi(\Tilde\xi,\cdot)\Big|_{[0,T]}\Big)\\[4mm]
&\ds\leq~\bigl\|\chi(\xi,\cdot) - \chi_n(\xi,\cdot)\|_{\L^\infty ([0,T])} + 
\delta\Big(\chi_n(\xi,\cdot)\Big|_{[0,T]}, \chi_n(\Tilde\xi,\cdot)\Big|_{[0,T]}\Big)
\\[4mm]
& \qquad \ds +\bigl\|\chi(\Tilde\xi,\cdot) - 
\chi_n(\Tilde \xi,\cdot)\|_{\L^\infty ([0,T])} \\[3mm]
&\ds\leq~\delta\Big(\chi_n(\xi,\cdot)\Big|_{[0,T]}, \chi_n(\Tilde\xi,\cdot)\Big|_{[0,T]}\Big)+{2\sigma\over 3}\,. \enda\eeq
This proves that the two paths $\chi_n(\xi,\cdot)$ and 
$\chi_n(\Tilde\xi,\cdot)$, which are parameterized by arc-length,
 cannot coincide over the entire interval $[0,T]$.
\v
{\bf 12.}  We are finally ready to prove (\ref{lsc2}). Let $\ve>0$ be given.
Since the weights $\Tilde W_i$ are uniformly bounded, by choosing $h>0$ small enough 
for every $i=1,\ldots,N$ we achieve
\bel{0045} \int_{a_i}^{b_i} \psi(\Tilde W_i(s))\, ds-\sum_{k=1}^{n(i)} \int_{s_{k-1}+h}^{s_k - h}\psi(\Tilde W_i(s))\, ds ~<~{\ve\over N}\,. \eeq
Since $\ve>0$ is arbitrary, to prove (\ref{lsc2}), it is thus suffices to show that
\bel{0047} \sum_{i = 1}^N\sum_{k=1}^{n(i)}\int_{s_{k-1}+h}^{s_k - h}\psi(\Tilde W_i(s))\, ds~\leq~\liminf_{n\to \infty} \E^{W,\psi}(\chi_n). \eeq
As shown in step {\bf 9}, there is a map
\bel{gi2}
\gamma_{i,k}~\mapsto~\Gamma_{i,k}~\doteq~\bigl\{\gamma^\sharp_l\,;~l=l(i,j,k)\bigr\},\eeq
which associates to the portion of elementary path $\gamma_{i,k}$ of $\chi$ 
a corresponding family  of $\ve_n$-good paths of $\chi_n$, as in Fig.~\ref{f:ir94}.
Using the ordering (\ref{ikprec}), by induction on $(i,k)$ we will show that
\bel{mest}
\int_{s_{k-1}+h}^{s_k - h}\psi(\Tilde W_i(s))\, ds ~\leq~
\sum_{\gamma_l^\sharp\in \Gamma_{i, k}} \int_{a_l^\sharp}^{b_l^\sharp} \psi(W_l^\sharp(s))\, ds,
\eeq
for every $(i,k)$.   By showing that paths $\gamma^\sharp_l$ belonging to distinct 
families $\Gamma_{i, k}$, $\Gamma_{i^\dagger, k^\dagger}$ are disjoint, we will conclude
\bel{sc6} \sum_{i = 1}^N\sum_{k=1}^{n(i)}\int_{s_{k-1}+h}^{s_k - h}\psi(\Tilde W_i(s))\, ds~\leq~
\sum_{i = 1}^N\sum_{k=1}^{n(i)}\sum_{\gamma_l^\sharp\in \Gamma_{i, k}} \int_{a_l^\sharp}^{b_l^\sharp} \psi(W_l^\sharp(s))\, ds~\leq~ \E^{W,\psi}(\chi_n),\eeq
for every $n\geq 1$ sufficiently large.  This will prove (\ref{0047}), and hence (\ref{lsc2}).
\v
{\bf 13.} In this step we prove our claim that paths $\gamma^\sharp_l$ belonging to distinct 
families $\Gamma_{i, k}$, $\Gamma_{i^\dagger, k^\dagger}$ are disjoint.
Assume $(i,k)\not= (i^\dagger, k^\dagger)$.   Two cases can occur.
\v
CASE 1: the elementary paths $\gamma_i, \gamma_{i^\dagger}$ are not contained in the same maximal
$\ve_0$-good path of~$\chi$.

In this case, there exists two distinct maximal paths $\Hat \gamma_p, \Hat \gamma_q$, which bifurcate at time
\bel{3001}\tau_{pq}~\leq~\min\{ s_{k-1}, \, s_{k^\dagger -1}\}\eeq
and such that
\bel{1016} \gamma_i(t)~=~\Hat\gamma_p(t),\quad\forall t\in [s_{k-1},s_k],\qquad 
\gamma_{i^\dagger}(t)~=~\Hat\gamma_q(t),\quad\forall t\in [s_{k^\dagger-1}, s_{k^\dagger}].  \eeq
If now $\xi\in A_{i}^j(s_k)$ and $\Tilde \xi \in A_{i^\dagger}^{j^\dagger}(s_{k^\dagger})$, 
by (\ref{+44}) there exists two $\ve_n$-good maximal paths for $\chi_n$ such that
\bel{3000}\bega{c}  \ds\chi_n(\xi,t)~=~\Hat\gamma'_j(t),\qquad \forall 0\leq t\leq \eta(\xi, s_k - h) ,     \\[3mm]
 \ds \chi_n(\Tilde\xi,t )~=~\Hat\gamma_{j^\dagger}'(t),\qquad \forall 0\leq t\leq \eta(\Tilde\xi, s_{k^\dagger}-h).  \enda\eeq
By (\ref{3001}) and the analysis in step {\bf 11}, the two paths $\chi_n(\xi,\cdot)$ and  $\chi_n(\Tilde\xi,\cdot)$ must bifurcate before time $\eta(\xi,s_{k-1} + h)\wedge \eta(\Tilde\xi, s_{k^\dagger-1} +h)$. 
Here and in the sequel we use the notation $a\wedge b\doteq\min\{a,b\}$.  
Calling $\tau'_{j j^\dagger}$ the time  where the two maximal paths $\Hat\gamma_j'$ and $\Hat\gamma'_{j^\dagger}$ bifurcate, by (\ref{3000}) and (\ref{+11}) one obtains
\bel{tjj} \tau'_{j j^\dagger}~\leq~ \eta(\xi,s_{k-1} + h)\wedge \eta(\Tilde\xi, s_{k^\dagger-1} +h)~\leq
~ s^{i,j}_{k-}\wedge s^{i^\dagger,j^\dagger}_{k^\dagger -}\,.\eeq

If now $A_i^j(s_k)$ is nonempty, by construction the path $\gamma^\sharp_l\in \Gamma_{i, k}$, which 
is contained in the maximal path $\Hat\gamma_j'$, is defined for $s\in [ s^{i,j}_{k-}, s_{k+}^{i,j} ]$. 
Similarly, 
$A_{i^\dagger}^{j^\dagger}(s_{k^\dagger})$ is nonempty, the path 
$\gamma^\sharp_{l^\dagger}\in \Gamma_{i^\dagger, k^\dagger}$ which is contained in 
$\Hat\gamma_{j^\dagger}'$ will be defined for 
$s\in [ s^{i^\dagger,j^\dagger}_{k^\dagger-}, s_{k^\dagger+}^{i^\dagger,j^\dagger} ]$. 

Since the two maximal $\ve_n$-good paths $\Hat\gamma_j'$ and $\Hat\gamma_{j^\dagger}'$
already bifurcate at the time (\ref{tjj}), the two paths $\gamma^\sharp_l$ and $\gamma^\sharp_{l^\dagger}$ are disjoint.
By the above argument, we conclude that
 the two families $\Gamma_{i, k}$ and $\Gamma_{i^\dagger,k^\dagger}$ consist of distinct paths. 
 \v
CASE 2: with  the partial ordering (\ref{ikprec}) one has $(i,k)\preceq(i^\dagger, k^\dagger)$. 

This implies that  there exists a maximal $\ve_0$-good path $\Hat\gamma_j$ in $\chi$, such that 
\bel{1011} \gamma_i(t)~=~\Hat\gamma_j(t),\quad\forall t\in [s_{k-1}, s_{k}],\qquad \gamma_{i^\dagger}(t)~=~\Hat\gamma_j(t),\quad\forall t\in [s_{k^\dagger-1}, s_{k^\dagger}] , \eeq
and moreover $s_k<s_{k^\dagger}$.
 For each fixed maximal $\ve_n$-good path $\Hat\gamma_j'$ in $\chi_n$, there are two cases:
\begin{itemize}
		\item $A_{i^\dagger}^j(s_{k^\dagger}) $ is nonempty. 
By (\ref{+44}) and (\ref{1011}) we thus have
	$A_{i^\dagger}^j(s_{k^\dagger})\subseteq A_i^j(s_k)$. Hence $A_i^j(s_k)$ is  nonempty as well.
By the definition (\ref{+11}) one has
	\bel{1015}\bega{rl}\ds s^{i^\dagger,j}_{k^\dagger-} - s^{i,j}_{k+}&\ds=~\inf_{\xi\in A_{i^\dagger}^j(s_{k^\dagger})} \eta(\xi,s_{k^\dagger-1}+h) - \inf_{\xi\in A_i^j(s_k)} \eta(\xi, s_k - h)\\[3mm]
	&\ds\geq \inf_{\xi\in A_{i^\dagger}^j(s_{k^\dagger})} \left[  \eta(\xi,s_{k^\dagger-1}+h) - \eta(\xi, s_k - h) \right]\\[3mm] 
	&\ds\geq~\inf_{\xi\in A_{i^\dagger}^j(s_{k^\dagger})}\left[  \eta(\xi,s_{k}+h) - \eta(\xi, s_k - h) \right]~\geq~2h.  \enda\eeq
	In this case, for every path $s\mapsto \gamma_{l^\dagger}^\sharp(s)$ 
	in $ \Gamma_{i^\dagger, k^\dagger}$, which is contained in $\Hat\gamma'_j$, the arc-length parameter ranges in $ [s_{k^\dagger-}^{i^\dagger,j} , s_{k^\dagger+}^{i^\dagger,j}]  $. On the other hand,  for every  path $s\mapsto \gamma_{l}^\sharp(s)$ in  $ \Gamma_{i, k}$, which is contained in $\Hat\gamma'_j$, the time parameter ranges in $ [s_{k-}^{i,j}, s_{k+}^{i,j}]  $. 
By (\ref{1015}) these two paths are disjoint.
	
\item $A_{i^\dagger}^j(s_{k^\dagger})$ is  empty. By construction, this implies that every path $\gamma_{l^\dagger}^\sharp\in  \Gamma_{i^\dagger,k^\dagger}$ is not contained in the maximal path $\Hat\gamma_j'$.  Thus, if $\gamma_l^\sharp\in \Gamma_{i, k}$ is contained in $\Hat\gamma_j'$,   
$\gamma_l^\sharp$ is disjoint from all the paths in $\Gamma_{i^\dagger,k^\dagger}$.
\end{itemize} 
Since the above analysis applies to each maximal path $\Hat\gamma_j'$ in $\chi_n$,	
	  we conclude that when $(i,k)\preceq (i^\dagger,k^\dagger)$, the two families 
	  $\Gamma_{i, k}$ and $\Gamma_{i^\dagger,k^\dagger}$ consist of disjoint paths.
\v
{\bf 14.} As before, let $\gamma_1,\ldots,\gamma_N$ be the elementary $\ve_0$-good paths in 
$\chi$. The weights $\Tilde W_i$ are then constructed along each $\gamma_i$ by the same 
inductive procedure as in (\ref{W1})--(\ref{posti}), for $i\in \I_p$, $p=1,2,\ldots$
We recall that $\I_p$ are the sets of indices introduced at (\ref{i1n}).

Toward a proof of  (\ref{0047}) we claim that, 
for any $i$, $1\leq k\leq n(i)$,
\bel{1020} \sum_{\gamma_l^\sharp\in \Gamma_{i, k}} \int_{a_l^\sharp}^{b_l^\sharp} \psi(W_l^\sharp(s))\, ds~\geq~\int_{s_{k-1}+h}^{s_k - h} \psi(\Tilde W_i(s))\, ds, \eeq
\bel{1021} \sum_{\gamma_l^\sharp\in \Gamma_{i, k},~ a_l^\sharp= s_{k-}^{i,j}} W_l^\sharp (a_l^\sharp)~\geq~\Tilde W_i(s_{k-1}+h).    \eeq
The above inequalities will be proved first for $i\in \I_1$ (i.e., for the outer-most branches),
then inductively for $i\in \I_2,\I_3,\ldots$

We begin by considering 
an elementary path $\gamma_i$ with $i\in \I_1$.  We compare the weight $\Tilde W_i$
along $\gamma_i$ with the sum of weights along the corresponding $\ve_n$-good 
paths $\gamma^\sharp_l$ of 
$\chi_n$.
On the last subinterval $[ s_{n(i)-1}+h,\, s_{n(i)}-h ]$ of $\gamma_i$, 
by (\ref{0032})-(\ref{0031}) the assumptions in Lemma \ref{l:mpc} are satisfied.
From (\ref{0054})-(\ref{1004}) we thus have
\bel{2001} \sum_{\gamma_l^\sharp\in \Gamma_{i, n(i)}} \int_{ a_l^\sharp}^{ b_l^\sharp}\psi(W^\sharp_l(s))\, ds~\geq~\int_{s_{n(i)-1} + h}^{s_{n(i)}-h} \psi(\Tilde W_i(s))\, ds\,, \eeq
\bel{2002} \sum_{ \gamma_l^\sharp\in \Gamma_{i, n(i)},~  a_l^\sharp = s_{n(i)-}^{i,j} } W_l^\sharp(a_l^\sharp )~\geq~\Tilde W_i( s_{n(i)-1} + h  ) . \eeq
Now consider the previous interval  $[ s_{n(i) - 2} +h,\, s_{n(i)-1}-h  ]$. By (\ref{TMI})-(\ref{TMI1}) 
it follows
\bel{2003}  \Tilde W_i( s_{n(i)-1} - h )~=~\Tilde W_i( s_{n(i) - 1} + h ) + \Tilde m_i(s_{n(i) - 1} -h ) -\Tilde m_i(s_{n(i)-1} + h)  .  \eeq
Hence by (\ref{0031}) and (\ref{2002}) one has
\bel{2004}\bega{rl}\ds \sum_{\gamma_l^\sharp\in \Gamma_{i,n(i)-1},~ b_l^\sharp = s^{i,j}_{(n(i)-1)+} } W_l^\sharp( b_l^\sharp )&\ds\geq~ \sum_{ \gamma_{l'}^\sharp\in \Gamma_{i, n(i)},~  a_{l'}^\sharp = s_{n(i)-}^{i,j} } W^\sharp_{l'}(a_{l'}^\sharp )\\[3mm] 
&\ds\qquad+~ \Tilde m_i( s_{n(i)-1} - h  ) -\Tilde m_i(s_{n(i) - 1} + h )\\[3mm]
&\ds \geq~ \Tilde W_i( s_{n(i)-1} - h ).\enda\eeq
By  (\ref{0032}) and (\ref{2004}) we can again apply Lemma \ref{l:mpc} on $\Gamma_{i,n(i)-1}$ and the restriction of $\gamma_i$ on $[s_{n(i) - 2} +h,\, s_{n(i)-1}-h  ]$. By similar arguments we prove (\ref{1020})-(\ref{1021}) for each $i\in \I_1, 1\leq k\leq n(i)$.
\v
{\bf 15.} Next, assume that (\ref{1020})-(\ref{1021}) 
have been proved for all $i\in \I_1\cup\cdots\cup\I_{p-1}$.
We claim that these same inequalities also hold for all $i\in \I_{p}$. 

Indeed, 
consider an elementary path $\gamma_i$ with $i\in \I_p$.  Along this path, consider the last 
subinterval, with  
$s\in [s_{n(i)-1}+h,\, s_{n(i)}-h]$.  
Recalling the construction of the weight $\Tilde W_i$ 
at (\ref{Wilast})-(\ref{posti}) and (\ref{TMI})-(\ref{TMI1}), we obtain
\bel{3002}\bega{rl}\ds \Tilde W_i(s_{n(i)} - h)&\ds=~\Tilde W_i(s_{n(i)})~=~\sum_{k\in \O(i)} \Tilde W_k(a_k+) + \left[\Tilde m_i(s_{n(i)}) - \sum_{k\in \O(i)} \Tilde m_k(a_k+) \right]\\[3mm]
&\ds =~\sum_{k\in \O(i)} \Tilde W_k(a_k + h) + \left[ \Tilde m_i(s_{n(i)} - h) - \sum_{k\in \O(i)} \Tilde m_k(a_k + h) \right].  \enda\eeq
\bel{3006}\bega{rl}\ds \sum_{\gamma_l^\sharp \in \Gamma_{i,n(i)},~b^\sharp_l = s_{n(i)+}^{i,j}  } W^\sharp_l(b_l^\sharp )&\ds\geq~ \sum_{k\in \O(i)} \sum_{\gamma^\sharp_{l'}\in \Gamma_{k,1},~a^\sharp_{l'}= s^{k,j}_{1-} } W^\sharp_{l'}(a^\sharp_{l'})\\[3mm]
&\ds\qquad + \left[ \Tilde m_i(s_{n(i)} - h)  - \sum_{k\in \O(i)} \Tilde m_k(a_k+h)   \right].  \enda\eeq
For each $k\in \O(i)$, the inductive assumption (\ref{2002}) yields
\bel{3003} \sum_{\gamma_{l'}^\sharp \in \Gamma_{k, 1},~ a_{l'}^\sharp = s^{k,j}_{1-}} W_{l'}^\sharp (a_{l'}^\sharp)~\geq~\Tilde W_k(a_k+h) .   \eeq
Hence, by (\ref{3002})--(\ref{3003}),
\bel{3004}  \sum_{\gamma^\sharp_l\in \Gamma_{i,n(i)},~b^\sharp_l = s^{i,j}_{n(i)+} } W^\sharp_l(b_l^\sharp)~\geq~\Tilde W_i(s_{n(i)} - h).   \eeq
Thanks to (\ref{3004})  and (\ref{0032}), we can use again  Lemma \ref{l:mpc} and conclude 
(\ref{2001})-(\ref{2002}). 
By backward induction on $k= n(i), n(i)-1,\ldots,1$, we then achieve the proof of (\ref{1020})-(\ref{1021})
as in step~{\bf 14}.   
\v
{\bf 16.} By induction on $\I_p$, $p=1,2,\ldots$, we conclude that the inequalities (\ref{1020}) 
hold for every $i=1,\ldots,N$ and every $k=1,\ldots, n(i)$. 
In turn, since  the families of paths $\Gamma_{i, k}$ are all disjoint from each other, from (\ref{1020})
we obtain (\ref{sc6}). As remarked in step {\bf 12}, this  implies the lower semicontinuity of the weighted irrigation cost.
\endproof

\section{Weights depending on the inclination of the branches}
\label{s:9}
\setcounter{equation}{0}
Aim of this section is to 
extend the previous results  to the case where the right hand side of the ODE 
in (\ref{ODE}) also depends on the inclination of the branch.
More precisely, if $s\mapsto\gamma(s)$ is a parameterization of the branch,
we replace (\ref{ODE}) with
\bel{fnew}W'(s)~=~-f\bigl(\dot\gamma(s), W(s)\bigr).\eeq
Concerning the function $f:\R^d\times\R\mapsto\R$, we shall assume
\begi\item[{\bf (A3)}] {\it The function $f=f(v,W)$ is continuous w.r.t.~both variables.
For each $v\in \R^d$, the map
$W\mapsto f(v,W)$ satisfies the same conditions as in (\ref{fp1}), namely
\bel{fp7}f(v,0)\,=\,0,\qquad  f_{W}(v,W)\,>\,0,
\qquad f_{WW}(v,W)\,\leq\,0\qquad \forall W>0.\eeq
For each $W> 0$, the map $v\mapsto f(v,W)$ is convex and positively homogeneous, namely
\bel{fp8}
f(r v, W)\,=\, r \, f(v,W)\qquad\forall r\geq 0.\eeq
}\endi

An example of a  function satisfying {\bf (A3)} is
$$f(v,W)~=~\left(|v_1| + \sqrt{v_1^2+v_2^2}\right) W^\beta,$$
where $0<\beta\leq 1$. 

Let now $\chi:[0,M]\times \R_+\mapsto \R^d$ be an irrigation plan.
When $f$ also depends on $v$, the weight functions 
$W(\xi,t)$ can be constructed following 
exactly the same procedure described in Section~\ref{s:23}.
The only difference is that, for each elementary path $\gamma_i:[a_i, b_i]\mapsto\R^d$,
the formulas    (\ref{W1})-(\ref{Wilast}) are now replaced respectively by
\bel{W15}
w(t)~=~\int_t^{b_i} f\bigl(\dot\gamma_i(s), w(s)\bigr)\, ds + m_i(t),\qquad\qquad t\in \,]a_i, b_i],
\eeq
\bel{Wil}
w(t)~=~\int_t^{b_i} f\bigl(\dot\gamma_i(s), w(s)\bigr)\, ds + m_i(t) + \ov w_i\,,\qquad\qquad t\in \,]a_i, b_i].\eeq
Here the upper dot denotes a derivative w.r.t.~the parameter $s$ along the arc.
We observe that all  conclusions of Lemma~\ref{l:mvode} remain valid if (\ref{ie3}) is replaced by
\bel{ie9} w(t)~=~\int_t^T f\bigl(s, w(s)\bigr)\, ds + m(t)\qquad\qquad\forall ~t\in [0,T],\eeq
assuming that $f:[0,T]\times [\ve, +\infty[$ is measurable w.r.t.~$t$ and Lipschitz continuous w.r.t.~$w$.

Relying on the assumptions {\bf (A3)}, the lower semicontinuity of the weighted irrigation cost 
proved in Theorem~\ref{t:51} can now be extended to this more general case.

\begin{theorem}\label{t:6} Consider a sequence $(\chi_n)_{n\geq 1}$ of irrigation plans,
all satisfying the assumption {\bf (A2)}, 
pointwise converging to an irrigation plan $\chi$. 
Assume that the function $\psi$ satisfies the conditions in {\bf (A1)}, while
$f$ satisfies {\bf (A3)}.
Then the corresponding
weighted costs satisfy
\bel{lsc4}
 \E^{W,\psi}(\chi)~\leq~\liminf_{n\to \infty} \E^{W,\psi}(\chi_n).\eeq
\end{theorem}

{\bf Proof.}    We shall follow step by step all the arguments in the proof 
 of Theorem~\ref{t:51}, and indicate only the modifications which are needed to cover
 this more general case.
 
 Steps {\bf 1--3}  and  {\bf 5--13} do not make any reference to the function $f$, and thus 
 remain valid without any change.   
 
 In step~{\bf 4} we considered an approximate family of weights $\Tilde W_i$ yielding almost the same cost
 as the original ones.  That construction must here be somewhat refined, approximating all 
 $\ve_0$-good paths in $\chi$ with polygonal lines.  That step is now replaced by
 \v
 {\bf 4$^\prime$.}
By the properties of $f$, there exist constants $L, \kappa$ such that 
\bel{fp4}
|f(v, w_1)- f(v, w_2)|~\leq~L\,|w_1-w_2|,\qquad\qquad\forall  w_1,w_2>\ve_0\,,\quad |v|\leq 1,\eeq
\bel{fp5}
|f(v,w)|~\leq~\kappa\qquad\qquad \forall  w\leq w_{max},\quad |v|\leq 1,\eeq
Here $w_{max}$ denotes the maximum weight over all $\ve_0$-good paths of $\chi$.

Let $\gamma_1,\ldots, \gamma_N$ be the elementary $\ve_0$-good paths in the irrigation plan $\chi$,
determined by the Path Splitting Algorithm {\bf (PSA)}, and let   $\ve>0$  be given.
By choosing $\delta,\delta_0>0$ sufficiently small, the following holds.

For each $\gamma_i:[a_i, b_i]\mapsto \R^d$, consider any set of intermediate times $s_k$ as in (\ref{+0})
with $$s_k-s_{k-1}~<~\delta_0\qquad\forall k\in \{1,\ldots, n(i)\}\,.$$
Define the piecewise constant multiplicity function
$\Tilde m_i(t)$ as in (\ref{TMI}).
Next, let $J_i\subset [a_i, b_i]$ by any measurable subset  with $meas(J_i)\leq \delta_0$ and define
\bel{TMI7} \tilde f_i(t,\omega)~\doteq~\left\{\bega{cl} f(\dot\gamma_i(t),\omega) - \delta_0 &\qquad 
\hbox{if}~~~t\notin J_i,\\[3mm]
0&\qquad \hbox{if} ~~~t\in J_i\,. \enda\right. \eeq
Then 
the corresponding weights $\Tilde W_i$ still satisfy (\ref{WTW}).

In the present case, an additional approximation will be useful.  Namely, we refine the partition (\ref{+0}) 
of the interval
$[a_i, b_i]$ by inserting points
\bel{part2}a_i\,=\,\tau_0\,<\,\tau_1\,<~\cdots~<\,\tau_{m(i)}\,=\,b_i\,,\eeq
and replace $\gamma_i$ with a path $\gamma_i^\diams$ which is affine on each sub-interval  $[\tau_{j-1}, \tau_j]$
and satisfies
$$\gamma_i(\tau_j) ~=~\gamma^\diams(\tau_j)\qquad\forall j.$$
Then we choose $h>0$ small enough and set 
$$J_i~\doteq~ \bigcup_j~ [\tau_j-h, \tau_j+h].$$
By choosing the partition (\ref{part2}) sufficiently fine, and $h>0$ sufficiently small, we can achieve
\bel{kl} \kappa L\cdot \sup_j |\tau_j-\tau_{j-1}|~<~{\delta_0\over 2}\,.\eeq
Setting
\bel{mds}m_i^\diams(s)~\doteq~m_i(s)-\delta_0\,,\eeq
the weights $W^\diams$ now satisfy
\bel{wds}
{d\over dt} \bigl[W^\diams(t) - m^\diams(t)\big]~=~- \tilde f(\dot \gamma^\diams(t), W^\diams(t)).
\eeq
If $\delta_0, h>0$ are chosen small enough, then 
the corresponding weight functions $W^\diams_i:[a_i, b_i]\mapsto\R_+$ satisfy
 the analogue of (\ref{WTW}), namely
 \bel{Wd}
 \|W_i^{\ve_0}- W^\diams_i\|_{\L^1([a_i,b_i])}~<~\ve,\qquad \qquad
|W_i^{\ve_0}(a_i+)-W^\diams_i(a_i+)|\,<\,\ve \,.\eeq

\begin{figure}[ht]
\c{\includegraphics[scale=0.5]{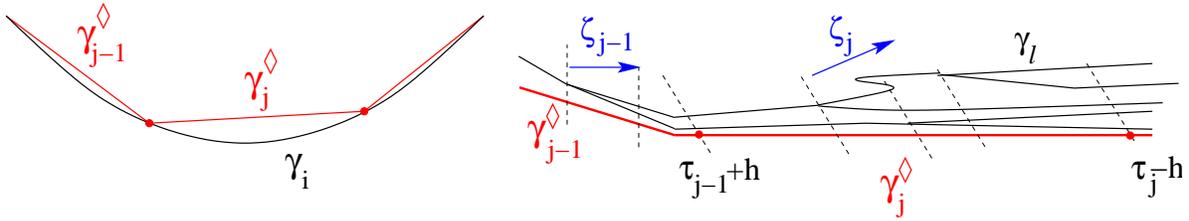}}
\caption{\small Left: approximating an $\ve_0$-good path $\gamma_{i}$ in the irrigation plan $\chi$
with a polygonal $\gamma^\diams$.  Right: the weight $W^\diams$ along the segment with 
endpoints $\gamma^\diams(\tau_{j-1}+h')$,  $\gamma^\diams(\tau_j-h')$ is compared with the sum of weights
along corresponding $\ve_n$-good paths $\gamma_l^\sharp$ of $\chi_n$.   Differently from 
the case illustrated in Fig.~\ref{f:ir94}, we now compare weights at points having the same inner product with
the vector $\zeta_j$, introduced at (\ref{904}).}
\label{f:ir99}
\end{figure}

Toward a proof of Theorem~\ref{t:6},
the heart of the matter is to achieve the inequalities (\ref{1020})-(\ref{1021}) in steps {\bf 14-15}.  
Since Lemma~\ref{l:mpc}
no longer applies, a different argument must now be developed.     The last three steps {\bf 14--16}
in the proof of Theorem~\ref{t:51} can be replaced by the steps below.
\v
{\bf 14$^\prime$.}
We wish to  compare
the weights $W_i^\diams$  with a sum of the weights along the corresponding $\ve_n$-good elementary paths
in the approximating irrigation plans $\chi_n$.  This will be done separately on each
 subinterval $[\tau_{j-1}, \tau_j]$, where  the tangent vector $\dot \gamma^\diams (t) =  v_j$ is constant.
 
As in step {\bf 9} of the previous proof,  in connection with $\gamma_i\Big|_{[\tau_{j-1}, \tau_j]}$ we can determine 
a family    $\Gamma_{i,j}$ of $\ve_n$-good paths 
$\gamma_\ell :[a_\ell , b_\ell ]\mapsto\R^d$ in the irrigation plan $\chi_n$, 
which approach $\gamma_i$ as $n\to\infty$ (see Fig.~\ref{f:ir99}, right).     
By construction, the corresponding weight and multiplicity functions 
$W_\ell, m_\ell :[a_\ell , b_\ell ]\mapsto\R$ are 
non-increasing and satisfy
\bel{dswm}{d\over ds} \bigl[W _\ell(s)- m_\ell(s)\bigr]~=~- f(\dot \gamma_\ell (s), W_\ell (s)).\eeq

 As remarked in the previous sections,  each elementary path $\gamma_\ell $ of  $\chi_n$ is contained in some 
 maximal path $\Hat \gamma_q:[a_q , b_q ]\mapsto\R^d$.   
 The set of elementary paths is partially ordered by setting
 $$\gamma_\ell ~\prec~\gamma_r $$
 if $\gamma_\ell $ and $\gamma_r $ are contained in the same maximal path, and 
 $b _\ell\leq a _r$.   The set of paths which bifurcate from the tip 
 of $\gamma_\ell$ is defined as
 $$\O(\ell)~\doteq~\{r\,;~~\gamma_\ell ~\prec~\gamma_r  ~~\hbox{and}~~
 b _\ell= a _r\}.$$
 At the endpoint of $\gamma_\ell$, by construction we have
 \bel{wbif}
 W_\ell(b_\ell)- m_\ell(b_\ell)~=~\sum_{r\in \O(\ell)} \bigl[W_r(a_r)- m_r (a_r)\bigr].\eeq
 Of course, the above definition here implies $a_r= b_\ell$.
 \v
 {\bf 15$^\prime$.}  Throughout this step we fix an elementary path $\gamma_i$ of the limit irrigation plan $\chi$.
 To shorten notation, we shall thus drop the index $i$ and simply write $\gamma^\diams=\gamma^\diams_i$,
 $W^\diams=W^\diams_i$, etc.
 To establish a comparison we observe that,
by the convexity and positive homogeneity of the map $v\mapsto f\bigl(v, W^\diams(\tau_{j})\bigr)$, 
there exists a vector $\zeta_j\in\R^d$ such that
\bel{904} f(v_j,W^\diams (\tau_{j})) ~=~\langle v_j, \zeta_j\rangle,\qquad\qquad f(v,W^\diams (\tau_{j}))~
\geq~\langle v, \zeta_j \rangle\quad \forall v\in \R^d\, .  \eeq
Notice that, if $f$ is smooth, then $\zeta_j=\nabla_v f\bigl(v_j, W^\diams(\tau_{j})\bigr)$ is simply the gradient 
of $f$ w.r.t.~the variable $v\in \R^d$. If $f$ is convex but not smooth, then $\zeta_j$ can be any sub-gradient.

By construction, for $t\in [\tau_{j-1},\tau_j]$ the derivative $\dot \gamma^\diams(t)=v_j$ is constant.
From (\ref{wds}), using   (\ref{TMI7}) and then (\ref{fp4})-(\ref{fp5}) and (\ref{kl}),  one obtains   
 \bel{ddw}\bega{l}\ds{d\over dt}
 \bigl[W^\diams (t)-m^\diams(t)\bigr]~\ds =~-\tilde f(v_j, W^\diams (t))\\[4mm]
 \qquad\ds \geq~-f\bigl(v_j, W^\diams (\tau_{j})\bigr) + \delta_0
 - \kappa L |\tau_j-t|
 ~\geq~-\langle v_j, \zeta_j\rangle + {\delta_0\over 2}\,.\enda\eeq

Let $\gamma _\ell:[ a _\ell,\,  b _\ell]\mapsto \R^d$ be any  $\ve_n$-good path in $\chi_n$.
 For  $t\in [\tau_{j-1}+h,\tau_j-h]$ we define
 \bel{sell}s_\ell(t)~\doteq ~\inf\Big\{s\geq a _\ell\,;~~\langle \zeta_j, \gamma_\ell (s)\rangle \,\geq\,
  \la \zeta_j, \gamma^\diams(t)
 \ra\Big\},\eeq
and consider the set of indices
$$I(t)~\doteq~\bigl\{\ell\,;~~a_\ell  \,<\, s_\ell(t)\,<\, b_\ell \bigr\}.$$
By an approximation argument, we can assume that each $\gamma_\ell$ is a polygonal, with 
$\langle\dot\gamma_\ell(t),\zeta_j\rangle\not= 0$ for a.e.~$t\in [a_\ell, b_\ell]$.  This implies
\bel{905}{d\over dt} s_\ell(t)~=~{\langle v_j, \zeta_j\rangle\over  \langle\dot\gamma_\ell(s_\ell(t)), \zeta_j\rangle}~>~0\eeq
for all except finitely many times $t$.  We can now estimate
 \bel{dwm}\bega{l}\ds{d\over dt}\sum_{\ell\in I(t)} \bigl[W_\ell(s_\ell(t))- m_\ell(s_\ell(t))\bigr]~=~ \ds-\sum_{\ell\in I(t)} 
 {\langle v_j, \zeta_j\rangle\over  \langle\dot\gamma_\ell(s_\ell(t)), \zeta_j\rangle}\cdot f(\dot\gamma_\ell(s_\ell(t)), W_\ell(s_\ell(t)))\\[4mm]
 \qquad \leq~\ds
 -\sum_{\ell\in I(t)}\langle v_j, \zeta_j\rangle\, {f(\dot\gamma_\ell(s_\ell(t)), W^\diams(\tau_{j})) \over 
 \langle\dot\gamma_\ell(s_\ell(t)), \zeta_j\rangle}\cdot \min\left\{{W_\ell(s_\ell(t))\over W^\diams(\tau_{j})}\,,\,1
 \right\}\\[4mm]
 \qquad \leq~\ds
 -\sum_{\ell\in I(t)} \langle v_j, \zeta_j\rangle\,\min\left\{{W_\ell(s_\ell(t))\over W^\diams(\tau_{j})}\,,\,1
 \right\}
\,.
 \enda
 \eeq
Here the first identity follows from (\ref{W15})-(\ref{Wil}) and (\ref{905}), 
 while the second inequality is a consequence of 
 (\ref{904}) and of the concavity of the 
map $W\mapsto f(v,W)$.   The third inequality follows from (\ref{904}).
Notice that,
 at points where two or more paths bifurcate, by (\ref{wbif}) the 
sum $\sum_{\ell\in I(t)}[W_\ell(t)-m_\ell(t)]$ remains continuous.  However, this sum will have
downward jumps at points where one of the maps $t\mapsto s_\ell(t)$ is discontinuous.   
\v
{\bf 16$^\prime$.} We are now ready to describe the  comparison argument that replaces
the estimates in step {\bf 15} of the proof of Theorem~\ref{t:51}.  
As in the previous proof, for each large $n$ we can identify a finite family  
of $\ve_n$-good paths
$\gamma_\ell$ in $\chi_n$ which converge to the elementary path $\gamma_i$ of $\chi$.
In particular, for $n\geq 1$ large enough, we can assume
\bel{mds}
\sum_\ell m_\ell(s_\ell(t))~\geq~m^\diams(t)\qquad\qquad \forall t\in [\tau_{j-1}+h\,,~\tau_j-h],~~j=1,2,\ldots, m(i).
\eeq
By backward induction, assume that
\bel{ww2}\bigl[W^\diams(\tau_j+h)-m^\diams(\tau_j+h)\bigr] ~\leq~\sum_{\ell\in I(\tau_j+h)} 
\bigl[ W_\ell(s_\ell(\tau_j+h))-m_\ell(s_\ell(\tau_j+h)\big].\eeq
Since $\tilde f=0$ for $t\in [\tau_j-h, \tau_j+h]$ while all differences $W_\ell-m_\ell$ are non-increasing,
this immediately yields
\bel{ww3}\bigl[W^\diams(\tau_j-h)-m^\diams(\tau_j-h)\bigr] ~=~\bigl[W^\diams(\tau_j+h)-m^\diams(\tau_j+h)\bigr] ~\leq~\sum_{\ell\in I(\tau_j-h)} 
\bigl[ W_\ell(s_\ell(\tau_j-h))-m_\ell(s_\ell(\tau_j-h)\big].\eeq
For $t\in [\tau_{j-1}+h\,,~\tau_j-h]$, we claim that the quantity
$$\Phi(t)~\doteq~\bigl[W^\diams (t)-m^\diams(t)\bigr]-\sum_{\ell\in I(t)} \bigl[W_\ell(s_\ell(t))- m_\ell(s_\ell(t))\bigr]$$
remains non-positive. Indeed,
 comparing the derivatives in (\ref{ddw}) and (\ref{dwm}) one obtains
\bel{dp1}\bega{l}\ds{d\over dt}\Phi(t)~\geq~\ds- \langle v_j, \zeta_j\rangle+{\delta_0\over 2} +\sum_{\ell\in I(t)} \langle v_j, \zeta_j\rangle\,\min\left\{{W_\ell(s_\ell(t))\over W^\diams(\tau_{j})}\,,\,1\right\}.
\enda
\eeq
Set
$$t^*~\doteq~\sup\,\bigl\{ t\leq \tau_j\,;~~\Phi(t)>0\bigr\}.$$
If $t^*>\tau_{j-1}+h$, to derive a contradiction we will show that 
\bel{dp2}
{d\Phi\over dt}\bigg|_{t=t^*}~>~0.\eeq
Toward this goal we observe that, by continuity, $\Phi(t^*)=0$. 
Since the map $t\mapsto W^\diams(t)$ is decreasing, by (\ref{mds}) we obtain
\bel{ww4}W^\diams(\tau_j) - \sum_\ell W_\ell(s_\ell(t^*))~\leq~W^\diams(t^*) - \sum_\ell W_\ell(s_\ell(t^*))~=~m^\diams(t^*) - \sum_\ell m_\ell(s_\ell(t^*))~\leq~0.\eeq
 By (\ref{dp1}) this implies
\bel{dp5} {d\Phi\over dt}\bigg|_{t=t^*}~\geq~\ds- \langle v_j, \zeta_j\rangle+{\delta_0\over 2} +\sum_{\ell\in I(t^*)} \langle v_j, \zeta_j\rangle\,\min\left\{{W_\ell(s_\ell(t^*))\over W^\diams(\tau_j)}\,,\,1\right\}~\geq~{\delta_0\over 2}.
\eeq
We thus conclude that $\Phi(t)\leq 0$ for all $t\geq\tau_{j-1}$.  This achieves the
key inductive step, showing that the inequality (\ref{ww2}) remains valid
with $j$ replaced by $j-1$.  

The remainder of the proof follows the same arguments used for Theorem~\ref{t:51}.
\endproof

\begin{remark}
{\rm  By a minor modification of the previous arguments,
the above results can be further extended 
to the case where $f=f(x,v,W)$  depends continuously also on the variable $x\in \R^d$.}
\end{remark}

\section{Optimal weighted irrigation plans}
\label{s:7}
\setcounter{equation}{0}
Given a positive, bounded Radon measure $\mu$ on $\R^d$, we define
$$\I^{W,\psi}(\mu)~\doteq~\inf_\chi \E^{W,\psi}(\chi).$$
where the infimum is taken over all irrigation plans for the measure $\mu$.
Relying on the lower semicontinuity of the weighted irrigation cost, proved in Theorems~\ref{t:51} and \ref{t:6},
we can now prove the existence of an optimal irrigation plan.

\begin{theorem}
	Let $\mu$ be a positive, bounded Radon measure on $\R^d$. Let 
$f$ satisfy the assumptions in {\bf (A3)} while $\psi$ satisfies {\bf (A1)}.   
 If $\mu$ admits an irrigation plan with bounded weighted cost, then there exists an irrigation plan with minimum weighted cost. \end{theorem}

{\bf Proof.} Let $M=\mu(\R^d)$ and let $(\chi_n)_{n\geq 1}$ be a minimizing sequence of irrigation plans
for $\mu$, so that
\bel{lew}\lim_{n\to \infty} ~\E^{W,\psi}(\chi_n)~=~\I^{W,\psi}(\mu).\eeq
Since $f\geq 0$, by construction the weights are larger than the corresponding 
multiplicities. Namely,
for every $\xi,t$ and $n\geq 1$ one has
$W_n(\xi,t)\geq m_n(\xi,t)$.
Since the costs in  (\ref{lew}) are bounded, we deduce
$$\int_0^M\int_0^{+\infty}|\dot\chi_n(\xi,t)|\, dt\, d\xi~\leq~C$$
for some constant $C$ and all $n\geq 1$.   

By the sequential compactness of traffic plans 
(see for example Proposition~3.27 in \cite{BCM}),
we can extract a subsequence $(\chi_{n_j})_{j\geq 1}$ 
pointwise converging to an irrigation plan $\chi$.  
The lower semicontinuity result proved in 
Theorem~\ref{t:6} yields 
\bel{lw7}\E^{W,\psi}(\chi)~\leq~\liminf_{n\to \infty} ~\E^{W,\psi}(\chi_n)~=~\I^{W,\psi}(\mu).\eeq
Hence $\chi$ achieves the minimum weighted cost.
\endproof
\v
We conclude by proving the lower semicontinuity of the weighted irrigation cost,
w.r.t.~weak convergence of measures.

\begin{theorem}\label{l:lsc}  Let 
$f$ satisfy the assumptions in {\bf (A3)} while $\psi$ satisfies {\bf (A1)}.   
	Let $(\mu_n)_{n\geq 1}$ be a sequence of bounded 
	positive measures,
	with uniformly bounded supports, weakly converging to $\mu$.  Then
	\bel{linf}\I^{W,\psi}(\mu)~\leq~\liminf_{n\to \infty} \I^{W,\psi}(\mu_n).\eeq
\end{theorem}
\v
{\bf Proof.} Without loss of generality, one can assume
\bel{1s}\liminf_{n\to \infty}~\I^{W,\psi}(\mu_n)~\doteq~K~<~+\infty\, . \eeq
	Let $\chi_n$ an optimal irrigation plan of $\mu_n$, so that
$\E^{W,\psi}(\chi_n)=\I^{W,\psi}(\mu_n)$ for every $n\geq 1$.
As in the previous proof, 
by sequential compactness we can extract a subsequence $(\chi_{n_j})_{j\geq 1}$ 
pointwise converging to an irrigation plan $\chi$.  A standard argument shows that 
$\chi$  provides an irrigation plan for the measure $\mu$. 
Using Theorem~\ref{t:51} we conclude
\bel{lw2}\I^{W,\psi}(\mu)~\leq~\E^{W,\psi}(\chi)~\leq~\liminf_{n\to \infty} ~\E^{W,\psi}(\chi_n)~=~\liminf_{n\to \infty} \I^{W,\psi}(\mu_n).\eeq
\endproof

\vs
{\bf Acknowledgment.}
This research was partially supported by NSF  grant DMS-1714237, ``Models of controlled biological growth".

\v

\end{document}